\documentclass[a4paper,12pt]{article}     

\usepackage{graphicx}		  
\usepackage{amsmath,amssymb}	

\newtheorem{lemma}{Lemma}[section]
\newtheorem{theorem}[lemma]{Theorem}
\newtheorem{proposition}[lemma]{Proposition}

\newtheorem{assumption}{Assumption}
\newtheorem{remark}[lemma]{Remark}
\newtheorem{fact}{Fact}
\newtheorem{note}[lemma]{Note}
\newtheorem{claim}{Claim}

\newcommand\LABEL[1]{\label{#1}}

\def\authorfont{\footnotesize}

\def\ccode#1{\par
\vspace*{8pt}
{\authorfont{\leftskip18pt\rightskip\leftskip
\noindent #1\par}}\par}

\newenvironment{Proof}{
\hspace*{-9mm}
{ \it Proof.}}
{\hfill {$\square$}\vspace{1.5em}}

\begin{document}

\begin{center}{
{\Large 
 The linear minimal 4-chart with three crossings}
\vspace{10pt}
\\ 
Teruo NAGASE and Akiko SHIMA
}
\end{center}

\begin{abstract}
Charts are oriented labeled graphs in a disk.
Every simple surface braid ($2$-dimensional braid) can be described by a chart.
Moreover, a chart represents an oriented closed surface
embedded in four-space.
In this paper, we study such embedded surfaces using charts.

Let $\Gamma$ be a chart.
Let ${\rm Cross}(\Gamma)$ denote the set of all the crossings of $\Gamma$, 
and let $\Gamma_m$ denote the union of all the edges of label $m$.
A 4-chart $\Gamma$  is said to be {\it linear}  
if every connected component of 
$$(\Gamma_1\cup \Gamma_3)-{\rm Cross}(\Gamma)$$
is acyclic. 
In this paper, we prove that 
after deleting free edges and hoops, 
every linear minimal $4$-chart with three crossings is 
lor-equivalent (equivalent under label changes, orientation reversals, and reflections) to
 the chart describing a $2$-twist-spun trefoil knot.
\end{abstract}
%
%
%

\ccode{2020 Mathematics Subject Classification. Primary 57K45,05C10; Secondary 57M15.}
\ccode{ {\it Key Words and Phrases}. surface-link, chart, C-move, crossing. }


\setcounter{section}{0}
\section{Introduction}


Charts are oriented labeled graphs in a disk (see  \cite{KnottedSurfaces},\cite{BraidBook}, and see Section~\ref{s:Prel}  for the precise definition of charts).
Let $D_1^2, D_2^2$ be 2-dimensional disks.
Any simple surface braid (2-dimensional braid) can be described 
by a chart,
here a simple surface braid is a properly embedded surface
$S$ in the 4-dimensional disk $D_1^2\times D_2^2$ such that
a natural map $\pi:S\subset D_1^2\times D_2^2\to D_2^2$ is 
a simple branched covering map of $D_2^2$ and
the boundary $\partial S$ is a trivial closed braid in
the solid torus $D_1^2\times \partial D_2^2$
(see \cite{BraidThree}, \cite[Chapter 14 and Chapter 18]{BraidBook}).
Moreover, from a chart, 
we can construct a simple closed surface braid in 4-space ${\Bbb R}^4$. This surface is an oriented closed surface 
embedded in ${\Bbb R}^4$.
On the other hand, any oriented closed surface 
embedded  in ${\Bbb R}^4$ is ambient isotopic to a simple
closed surface braid
 (see \cite{BraidThree},\cite[Chapter 23]{BraidBook}). 
A C-move 
is a local modification between two charts
in a disk (see Section~\ref{s:Prel} for C-moves).
A C-move between two charts induces 
an ambient isotopy between oriented closed surfaces 
corresponding to the two charts.
In this paper, we study oriented closed surfaces in 4-space
by using charts.

We will work in the PL category or smooth category. All submanifolds are assumed to be locally flat.

A {\it surface-link} is a closed surface embedded in 4-space ${\Bbb R}^4$. A {\it $2$-link} is a surface-link each of whose connected component is a 2-sphere. An orientable surface-link is called a {\it ribbon} surface-link if there exists an immersion of a 3-manifold $M$ into ${\Bbb R}^4$ sending the boundary of $M$ onto the surface-link such that each connected component of $M$ is a handlebody and its singularity consists of ribbon singularities, here a ribbon singularity is a disk in the image of $M$ whose pre-image consists of two disks; one of the two disks is a proper disk of $M$ and the other is a disk in the interior of $M$. In the words of charts, a ribbon surface-link is a surface-link corresponding to a {\it ribbon chart}, a chart C-move equivalent to a chart without white vertices \cite{BraidThree}. A chart is called a {\it $2$-link} chart if a surface-link corresponding to the chart is a 2-link.

Kamada showed that any 3-chart is a ribbon chart \cite{BraidThree}. 
Kamada's result was extended by Nagase and Hirota: Any 4-chart with at most one crossing is a ribbon chart \cite{NagaseHirota}. We showed that any $n$-chart with at most one crossing is a ribbon chart \cite{OneCrossing},\cite{MinimalChart}. We also showed that any 2-link chart with at most two crossings is a ribbon chart 
\cite{TwoCrossingI},\cite{TwoCrossingII}.
We investigated the structure of $c$-minimal charts with
two crossings 
(see Section~\ref{s:Prel} for the precise definition of a $c$-minimal chart),
and gave an enumeration of the charts with two crossings in
\cite{StI}, \cite{StII}.

\begin{figure}[h]
\begin{center}
\includegraphics{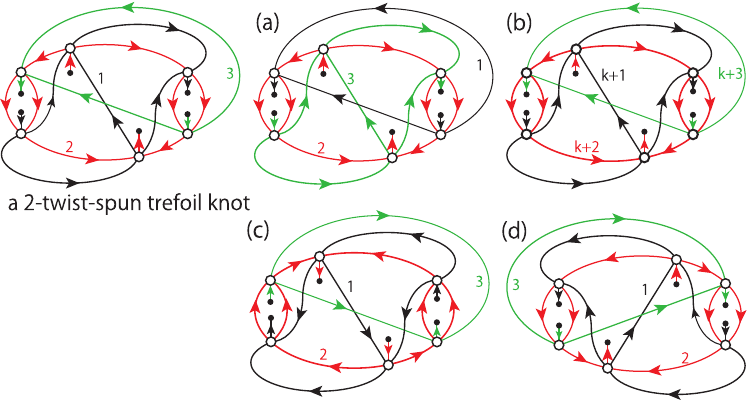}
\end{center}
\caption{\label{Fig01} Charts are lor-equivalent to the 4-chart describing a 2-twist-spun trefoil knot. Here $k$ is 
a positive integer.}
\end{figure}

Two charts are said to be 
{\it lor-equivalent}
(Label-Orientation-Reflection equivalent)
provided that
one of the charts is obtained from
the other by a finite sequence of 
the following five modifications:
\begin{enumerate}
\item[(i)]
consider an $n$-chart as an $(n+1)$-chart,
\item[(ii)]
for an $n$-chart,
change all the edges of label $k$ to 
ones of label $n-k$
for each $k=1,2,\cdots,n-1$,
simultaneously 
(see Fig.~\ref{Fig01}(a)),
\item[(iii)] 
add a positive constant integer $k$ 
to all the labels simultaneously 
(so that 
the $n$-chart changes to an $(n+k)$-chart)
(see Fig.~\ref{Fig01}(b)),
\item[(iv)]
reverse the orientation of 
all the edges (see Fig.~\ref{Fig01}(c)),
\item[(v)]
change a chart 
by the reflection in the sphere 
(see Fig.~\ref{Fig01}(d)).
\end{enumerate} 
For example, the charts as shown in 
Fig.~\ref{Fig01} are lor-equivalent to
the 4-chart describing a 2-twist-spun trefoil knot.

An edge in a chart is called 
a {\it free edge}
if it has
two black vertices.
A {\it hoop} is a closed edge of a chart without vertices
and hence also without crossings.

Let $\Gamma$ be a chart.
For each label $m$,  let $\Gamma_m$ denote
the union of all the edges of label $m$.

Let $\Gamma$ be a chart. Let ${\rm Cross}(\Gamma)$ denote the set of all crossings of the chart. Let $m$ be a label with 
$\Gamma_m \cap {\rm Cross}(\Gamma)\not= \emptyset$. 
Then $\Gamma_m$ is said to be {\it linear} provided that 
each connected component of $\Gamma_m - {\rm Cross}(\Gamma)$ is acyclic
(i.e. does not contain any cycle).

A chart $\Gamma$ is said to be {\it linear} provided that
any $\Gamma_m$ containing a crossing is linear.

We shall show the following theorem.

\begin{theorem} 
\label{MainTheorem} 
After deleting free edges and hoops,
any linear minimal $4$-chart with three crossings is 
lor-equivalent to the chart describing
a $2$-twist-spun trefoil knot.
\end{theorem}

The paper is organized as follows.
In Section~\ref{s:Prel},
we define charts and minimal charts.
In Section~\ref{s:mal-cycle},
we introduce a cycle consisting of four edges of label $m$, called a mal-cycle, which bounds a disk containing exactly one crossing and satisfying some condition. And we shall show that there does not exist a mal-cycle in any minimal chart.
In Section~\ref{s:ConsecutiveTripletLemma},
we review a useful lemma
called Consecutive Triplet Lemma (Lemma~\ref{ConsecutiveTripletLemma}).
In Section~\ref{s:AnacanthousBody},
we investigate 
the set  $AB(\Gamma_1)\cup AB(\Gamma_3)$ for a linear minimal 4-chart $\Gamma$,
where $AB(\Gamma_k)$ is obtained from $\Gamma_k$ by omitting terminal edges.
Here a terminal edge is an edge with a white vertex and a black vertex.
In Section~\ref{s:IO-Path}, we investigate simple arcs in $(AB(\Gamma_1)\cup AB(\Gamma_3))-{\rm Cross}(\Gamma)$.
In Section~\ref{s:4ChartTwoCrossing}, we shall show that 
after deleting free edges and
hoops,
any linear minimal 4-chart with two crossings
is lor-equivalent to the 4-chart describing a turned torus-link of a Hopf link. 
In Section~\ref{s:ProofTheorem},
we shall show the main theorem.



\section{Preliminaries}
\label{s:Prel}

In this section, 
we introduce 
the definitions of charts and related terminology.

\begin{figure}[htb]
\begin{center}
\includegraphics{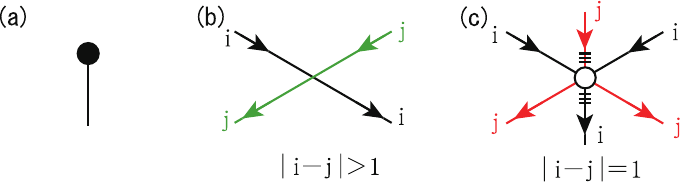}
\end{center}
\caption{ \label{Fig02} (a) A black vertex. (b) A crossing. (c) A white vertex. 
Each arc with three transversal short arcs is a middle arc at the white vertex. }
\end{figure}

Let $n$ be a positive integer.
An $n$-{\it chart}  
(a braid chart of degree $n$ \cite{KnottedSurfaces}
or a surface braid chart of degree $n$ \cite{BraidBook}) 
is 
an oriented labeled graph in the interior of a disk,
which may be empty 
or
have closed edges without vertices
satisfying the following four conditions
(see Fig.~\ref{Fig02}):
\begin{enumerate}
\item[(i)] 
Every vertex has degree $1$, $4$, or $6$.
\item[(ii)] 
The labels of edges are 
in $\{1,2,\dots,n-1\}$.
\item[(iii)]
In a small neighbourhood of
each vertex of degree $6$,
there are six short arcs,
three consecutive arcs are
oriented inward 
and
the other three are outward,
and
these six are labeled $i$ and $i+1$
alternately for some $i$,
where the orientation and label of
each arc are inherited from
the edge containing the arc.
\item[(iv)]
For each vertex of degree $4$,
diagonal edges have the same label
and
are oriented coherently,
and the labels $i$ and $j$ of
the diagonals satisfy $|i-j|>1$.
\end{enumerate}
We call a vertex of degree $1$ a {\it black vertex},
a vertex of degree $4$ a {\it crossing}, and 
a vertex of degree $6$ a {\it white vertex}
respectively.

Among six short arcs
in a small neighbourhood of
a white vertex,
a central arc of each three consecutive arcs
oriented inward (resp. outward) 
is called a   
{\it middle arc} at the white vertex
(see Fig.~\ref{Fig02}(c)).


\begin{figure}[hbt]
\begin{center}
\includegraphics{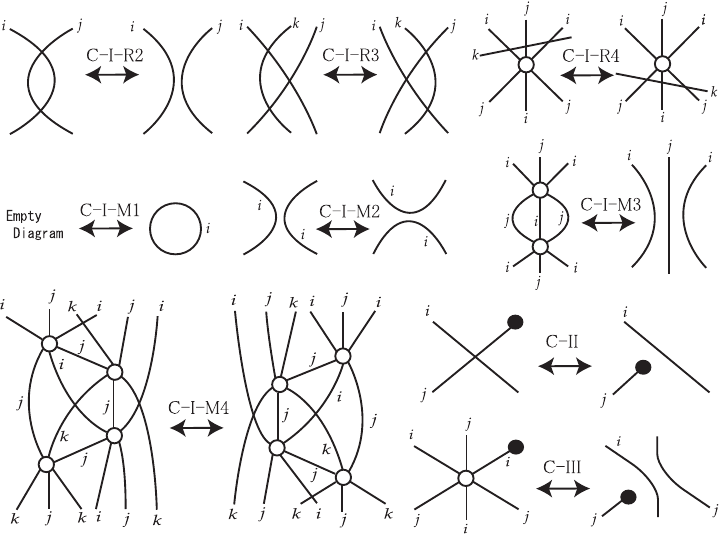}
\end{center}
\caption{ \label{Fig03} For the C-III move, 
the edge with the black vertex does not contain a middle arc at
a white vertex in the left figure. }
\end{figure}

Now {\it C-moves} are local modifications 
of charts as shown in Fig.~\ref{Fig03}
(cf. \cite{KnottedSurfaces}, 
\cite{BraidBook} and \cite{Tanaka}).
Two charts are said to be {\it C-move equivalent}  
if there exists a finite sequence of C-moves 
which modifies one of the two charts 
to the other.

\begin{figure}[bth]
\begin{center}
\includegraphics{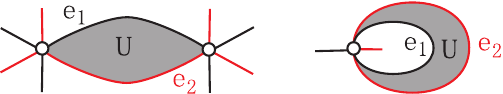}
\end{center}
\caption{ \label{Fig04} Bigons.}
\end{figure}

Let $\Gamma$ be a chart. 
Let $e_1$ and $e_2$ be edges of $\Gamma$
which connect two white vertices $w_1$ and $w_2$
where possibly $w_1=w_2$.
Suppose that 
the union $e_1\cup e_2$ bounds 
an open disk $U$.
Then the closure ${\rm Cl}(U)$ 
is called 
a {\it bigon} of $\Gamma$
provided that
any edge containing $w_1$ or $w_2$ 
does not intersect the open disk $U$
(see Fig.~\ref{Fig04}).
Note that neither $e_1$ nor $e_2$ contains a crossing.

%

Let $\Gamma$ be a chart.
Let 
$c(\Gamma),~w(\Gamma),~f(\Gamma)$, 
and 
$b(\Gamma)$ be 
the number of crossings, 
the number of white vertices, 
the number of free edges, 
and 
the number of bigons of $\Gamma$,
respectively.
The 4-tuple $(c(\Gamma),w(\Gamma),-f(\Gamma),-b(\Gamma))$ is called the 
{\it $c$-complexity} of the chart $\Gamma$
(see \cite{BraidThree} 
for complexities of charts).

A chart $\Gamma$ is said to be 
{\it $c$-minimal } if
its $c$-complexity is minimal among the charts 
which are C-move equivalent to 
the chart $\Gamma$
with respect to 
the lexicographic order of the 
4-tuple of integers.
If a chart is $c$-minimal, 
then we say that the chart is {\it minimal} in this paper.

A hoop is said to be {\it simple}
if one of the complementary domains of the hoop does not contain any white vertices.
An {\it oval nest} is a free edge 
together\index{oval nest} 
with some concentric simple hoops. 

An edge in a chart is called 
a {\it terminal edge}
if it has
a white vertex and a black vertex.

\begin{proposition}
\label{Assumption0}
{\rm (\cite[Remark~2.3]{MinimalChart}, 
\cite[Proposition~2.3]{StI})}
Let $\Gamma$ be a minimal chart in a disk $D^2$. 
Then we have 
the following:
\begin{enumerate}
\item[{\rm (a)}]
If an edge of $\Gamma$ contains a black vertex, 
then the edge is a terminal edge or a free edge.
\item[{\rm (b)}]
Any terminal edge of $\Gamma$ contains a middle arc 
at its white vertex.
\end{enumerate}
\end{proposition}

\begin{proposition}
\label{MoveOutSimpleHoop}
{\rm (\cite[Proposition~2.4]{StI})}
Let $\Gamma$ be a minimal chart in a disk $D^2$.
For any regular neighbourhood $N$ 
of $\partial D^2$ in $D^2$,
there exists a minimal chart $\Gamma'$ 
obtained from $\Gamma$ 
by C-I-M2 moves 
and ambient isotopies of $D^2$ such that
\begin{enumerate}
\item[{\rm (a)}]
$\Gamma'\cap (D^2-N)$ contains neither free edges
nor simple hoops,
\item[{\rm (b)}]
$\Gamma'\cap N$ consists of oval nests, 
simple hoops and free edges.
\end{enumerate}
\end{proposition}

In the following lemma,
we investigate the difference of a chart in a disk and in a 2-sphere.
This lemma follows from the existence of a natural one-to-one
correspondence between $\{$charts in $S^2\}/$C-moves and 
 $\{$charts in $D^2\}/$C-moves, conjugations (\cite[Chapter 23 and Chapter 25]{BraidBook}) 

\begin{lemma}
\label{MoveInfty}
{\rm (\cite[Lemma 2.1]{ChartApp1})}
Let $\Gamma$ and $\Gamma'$ be charts in a disk $D^2$.
Suppose that $\Gamma$ is ambient isotopic to $\Gamma'$ in the one-point compactification of the open disk ${\rm Int}(D^2)$,
i.e. the $2$-sphere.
Then there exist hoops $C_1,C_2,\dots,C_k$ in ${\rm Int}(D^2)$ such that 
\begin{enumerate}
\item[{\rm (a)}] the chart $\Gamma$ is obtained from $\displaystyle{\Gamma'\cup(\bigcup^k_{i=1}C_i)}$ by C-moves in the disk,
\item[{\rm (b)}] the chart $\Gamma'$ and hoops $C_1,C_2,\dots,C_k$ 
are mutually disjoint, and
\item[{\rm (c)}] each hoop $C_i$ bounds a disk containing the chart $\Gamma'$ in the disk $D^2$.
\end{enumerate}
\end{lemma}

To make the argument simple, we assume that 
the charts lie on the 2-sphere instead of the disk.
\begin{assumption}
In this paper,
all charts are contained in the $2$-sphere $S^2$.
\end{assumption}
We have the special point in the 2-sphere $S^2$, called the point at infinity,
 denoted by $\infty$.
In this paper, all charts are contained in a disk such that the disk 
does not contain the point at infinity $\infty$.

For any minimal chart in the $2$-sphere $S^2$
we can move free edges and simple hoops into 
a regular neighbourhood $U_\infty$ of the point at infinity 
$\infty$ in  $S^2$ 
by C-I-M2 moves and ambient isotopies of $S^2$ 
keeping fixed the point $\infty$
by Proposition~\ref{MoveOutSimpleHoop}.
During argument,
if free edges or simple hoops appear, 
we immediately move them 
into a regular neighbourhood $U_\infty$ of the point 
$\infty$ in $S^2$.
Thus we assume the following 
(cf. \cite{ChartApp1},\cite{OneCrossing}, \cite[Assumption 1]{MinimalChart}):

\begin{assumption}
\label{NoSimpleHoop}
All free edges and simple hoops in $\Gamma$ 
are moved into a small neighbourhood $U_\infty$ 
of the point at infinity $\infty$. 
Hence
we assume that 
$\Gamma$ does not contain free edges
nor simple hoops, 
unless otherwise stated. 
\end{assumption}

By Proposition~\ref{Assumption0} and Lemma~\ref{MoveInfty},
we assume the following:

\begin{assumption}
\label{AssumeTerminal}
If an edge of $\Gamma$
contains a black vertex,
then the edge is a free edge 
or a terminal edge.
Moreover 
any terminal edge contains a middle arc.
\end{assumption}

\begin{assumption}
\label{Infinity}
The point at infinity $\infty$ is moved in any complementary domain of $\Gamma$.
\end{assumption}

The following lemma will be used in
the proof of the main theorem. 

\begin{figure}
\begin{center}
\includegraphics{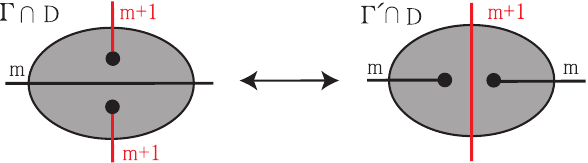}
\end{center}
\caption{ \label{Fig05} The gray regions are disks $D$, and $m$ is a label.}
\end{figure}

\begin{lemma}
{\rm (Cut Edge Lemma)(\cite[Lemma 18.24 (E)]{BraidThree},\cite[Lemma 5.1]{ChartApp1})}
\label{CutEdgeLemma}
Let $\Gamma$ and $\Gamma'$ be charts, and $D$ a disk with $\Gamma \cap D^c=\Gamma' \cap D^c$. 
If $\Gamma \cap D$ and $\Gamma' \cap D$ are sets as shown in Fig.~\ref{Fig05}, and if both $\Gamma_{m+1}\cap D$ and $\Gamma_{m}'\cap D$ consist of an inward arc and an outward arc,
then $\Gamma$ is C-move equivalent to $\Gamma'$. 
\end{lemma}

%

\begin{lemma}
{\rm (\cite{NagaseHirota})}
\label{4ChartRibbon}
Any $4$-chart with at most one crossing is a ribbon chart.
\end{lemma}

\begin{lemma}
\label{4ChartSimpleHoop}
Let $\Gamma$ be a minimal $4$-chart
with at most three crossings.
Then any hoop is a simple hoop.
\end{lemma}

\begin{Proof}
Let $C$ be a hoop.
Then $C$ separates $S^2$ into two disks,
say $D_1$ and $D_2$.
Let $\Gamma^*=\Gamma\cap {\rm Int}(D_1)$ and
$\Gamma^{**}=\Gamma\cap {\rm Int}(D_2)$.

Since $\Gamma$ contains at most three crossings,
one of  $\Gamma^*$ and $\Gamma^{**}$ contains at most one crossing.
Thus by Lemma~\ref{4ChartRibbon},
the subchart is a ribbon chart.
Therefore $C$ is a simple hoop.
\end{Proof}

In this paper
for a subset $X$ in a space
we denote 
the interior of $X$,
the boundary of $X$ and
the closure of $X$
by ${\rm Int}(X)$, $\partial X$
and ${\rm Cl}(X)$
respectively.



\section{Mal-cycles}
\label{s:mal-cycle}

In this section,
we introduce a cycle consisting of four edges of label $m$, called a mal-cycle, which bounds a disk containing exactly one crossing satisfying some condition. And we shall show that there does not exist a mal-cycle in any minimal chart.

Let $\Gamma$ be a chart and $e$ an edge with two vertices $v_1,v_2$.
If the edge $e$ is oriented from $v_1$ to $v_2$, we say that 
the edge $e$ is {\it outward at $v_1$} and also that 
the edge $e$ is {\it inward at $v_2$}.

Let $\Gamma$ be a chart and $m$ a label of the chart.
Let $\ell$ be an arc in $\Gamma_m$ each of whose boundary points is 
a white vertex or a crossing.
If any vertex in the interior of the arc $\ell$ is contained 
in a terminal edge of label $m$, 
then we call $\ell$ an {\it anacanthous path} of label $m$.

Let $\Gamma$ be a chart, $D$ a disk, and 
$v$ a vertex on the boundary of the disk $D$. 
An edge $e$ is called  an {\it inner edge for the disk $D$ at the vertex $v$} provided that 
the closure of $e \cap {\rm Int}(D)$ contains the vertex $v$ here 
${\rm Int}(D)$ is the interior of the disk $D$. 
If the edge $e$ is oriented inward (resp. outward) at the vertex $v$, 
then the edge $e$ is called
  an {\it inner edge for the disk $D$ 
inward $($resp. outward$)$ at the vertex $v$}.

If a vertex is contained in the interior of an arc $L$, 
then the vertex is called an {\it inner vertex} of the arc $L$.

\begin{figure}[hbt]
\begin{center}
\includegraphics{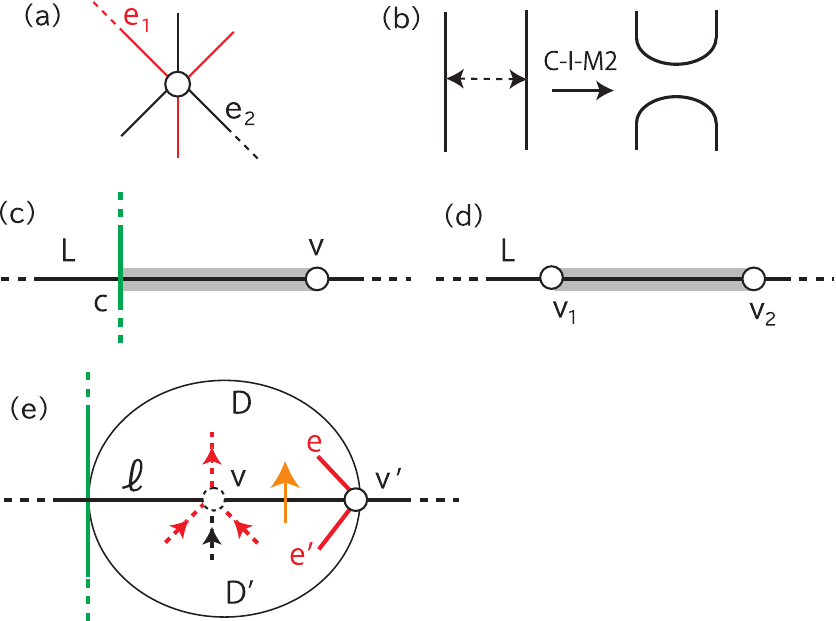}
\end{center}
\caption{ \label{Fig06} 
(a) Neither $e_1$ nor $e_2$ is a terminal edge.
(b) The double-headed dotted arrow indicates that 
a C-I-M2 move will be applied.
(c) There is no white vertex between the vertex $c$ and the vertex $v$ on the arc $L$.
(d) There is no white vertex between the vertex $v_1$ and the vertex $v_2$ on the arc $L$.
(e) A direction indicator for the arc $\ell$.}
\end{figure}

In constructing a chart in this paper, we use the following symbols 
in Fig.~\ref{Fig06} throughout this paper. 

\begin{enumerate}
\item[(1)] 
In arguments, if we do not decide an edge $e$ is a terminal edge or not, 
we draw just a short line segment for the edge $e$.
If we are confident of an edge $e$ not a terminal edge 
but still do not know the other vertex of the edge, 
then we add two dots or more at the end of the arc 
(see Fig.~\ref{Fig06}(a)).
\item[(2)]  We draw a double-headed dotted arrow between arcs which we will apply a C-I-M2 move (see Fig.~\ref{Fig06}(b)).

\item[(3)] If an arc does not contain any white vertex in its interior, 
we put a thick gray line under the arc 
(see Fig.~\ref{Fig06}(c),(d)). 

\item[(4)] Let $\ell$ be an anacanthous path and 
$D,D'$ disks with $D \cap D' =\ell$. 
To indicate that there does not exist any inner edge for $D$ inward at an inner vertex of the arc $\ell$
or that there does not exist any inner edge for $D'$ outward 
at an inner vertex of the arc $\ell$
we put a short arrow transverse to the arc $\ell$ 
pointing to the disk $D$ from the disk $D'$ (see Fig.~\ref{Fig06}(e)).
This arrow is called a {\it direction indicator for the arc $\ell$}. 
We also say that the direction indicator points to $D$ from $D'$.
In Fig.~\ref{Fig06}(e), 
if there exists an inner vertex $v$ of 
the arc $\ell$, 
then any inner edge of $D$ containing $v$ is outward at 
$v$ and any inner edge of $D'$ containing $v$ is inward at $v$.
But we do not know the two edges $e,e'$ are inward or outward at $v'$.
\end{enumerate}

%

\begin{figure}[hbt]
\begin{center}
\includegraphics{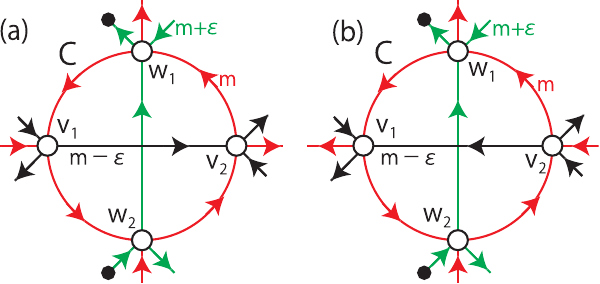}
\end{center}
\caption{ \label{Fig07} Mal-cycles.}
\end{figure}

Let $\Gamma$ be a chart and $m$ a label of the chart. Suppose $\varepsilon\in\{+1,-1\}$. 
Let $D$ be a disk whose boundary $C$ is a cycle of label $m$ 
with four vertices $v_1,w_1,v_2,w_2$ situated in this order 
on the cycle $C$.
The cycle $C$ is called a {\it mal-cycle} of label $m$ provided that
\begin{enumerate}
\item[(i)] the cycle $C$ is oriented clockwise or counterclockwise by the orientation inherited from the chart,
\item[(ii)]
 $v_1,v_2 \in \Gamma_m \cap \Gamma_{m-\varepsilon}$, and
$w_1,w_2 \in \Gamma_m \cap \Gamma_{m+\varepsilon}$,
\item[(iii)] each of $w_1,w_2$ is a vertex of a terminal edge 
of label $m+\varepsilon$ outside the disk $D$, and
\item[(iv)] the disk $D$ contains exactly one crossing 
and the four white vertices $v_1,v_2,w_1,w_2$ (see Fig.~\ref{Fig07}).
\end{enumerate}

\begin{lemma}
\label{MalCycleLemma}
There does not exist a mal-cycle in any minimal chart.
\end{lemma}

\begin{figure}[hbt]
\begin{center}
\includegraphics{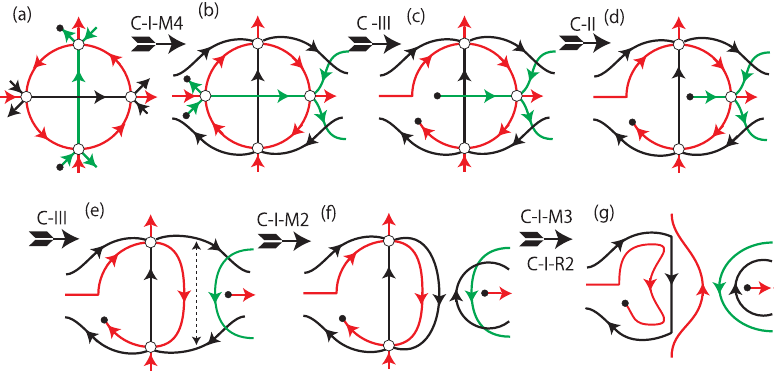}
\end{center}
\caption{ \label{Fig08} 
The double-headed dotted arrow in (e) indicates 
the site where a C-I-M2 move will be applied.}
\end{figure}

\begin{Proof}
Suppose that there exists a mal-cycle in a minimal chart.
We only show the case shown in Fig.~\ref{Fig07}(a).
Apply a C-I-M4 move to a mal-cycle 
(see Fig.~\ref{Fig08}(a),(b)). 
Then, as shown in Fig.~\ref{Fig08}, applying a C-III move, 
a C-II move, 
a C-III move, a C-I-M2 move, a C-I-M3 move and a C-I-R2 move one by one, we can eliminate 
one crossing and four white vertices.
This contradicts the minimality of the chart. 
\end{Proof}

\begin{figure}[hbt]
\begin{center}
\includegraphics{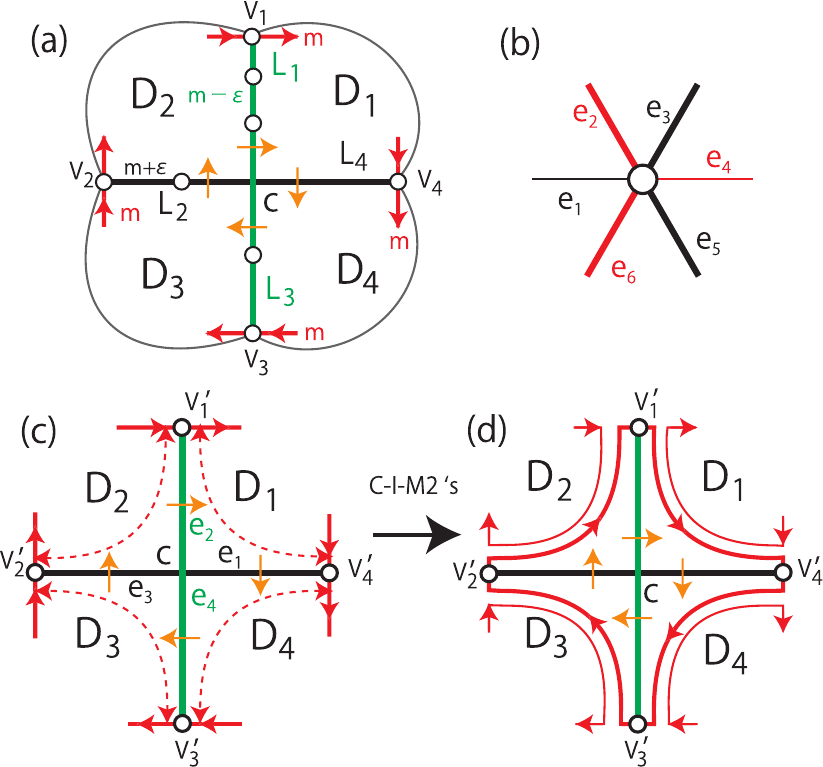}
\end{center}
\caption{ \label{Fig09} 
(a) A pinwheel.
(b) $e_2,e_6$ are side edges of $e_1$, and
$e_3,e_5$ are side edges of $e_4$.
(c) A neighbourhood of the crossing $c$.
(d) A mal-cycle.}
\end{figure}

Let $\Gamma$ be a chart, $m$ a label of the chart, 
and $c$ a crossing of the chart.
Let $\varepsilon\in\{+1,-1\}$. 
Let $L_1,L_3$ be anacanthous paths of label $m-\varepsilon$ and
$L_2,L_4$ anacanthous paths of label $m+\varepsilon$ such that 
\begin{enumerate}
\item[(i)]
for $i=1,2,3,4$ the arc $L_i$ contains the crossing $c$, 
and
a white vertex $v_i$,
\item[(ii)] for each different pair $L_i, L_j$ of $L_1,L_2,L_3,L_4$, 
we have $L_i\cap L_j=c$, and
\item[(iii)] each of $v_1$ and $v_3$ is contained in a terminal edge 
of label $m-\varepsilon$, or
each of $v_2$ and $v_4$ is contained in a terminal edge 
of label $m+\varepsilon$.
\end{enumerate}
The set $L_1 \cup L_2 \cup L_3 \cup L_4$ is called 
a {\it pinwheel} provided that
there exist four disks $D_1,D_2,D_3,D_4$ 
situated counterclockwise or clockwise
around the crossing $c$ in this order such that
\begin{enumerate}
\item[(i)]
$D_1\cap D_3 = D_2\cap D_4=c$ and 
$D_i\cap D_{i+1}=L_i$ for $i=1,2,3,4$ 
where we understand $D_5=D_1$,
\item[(ii)] $D= D_1\cup D_2 \cup D_3 \cup D_4$ does not contain 
any crossing except the crossing $c$, and
$\Gamma \cap D \subset \Gamma_{m-\varepsilon}\cup \Gamma_m\cup \Gamma_{m+\varepsilon}$,
\item[(iii)] 
\begin{enumerate}
\item[(a)]
for each $i=1,2,3,4$, the direction indicator for the arc $L_i$ pointing to $D_i$ 
from $D_{i+1}$ and there exist an inner edge of label $m$ 
for $D_i$ outward at $v_i$ and 
an inner edge of label $m$ for $D_{i+1}$ 
inward at $v_i$ or
\item[(b)] 
for each $i=1,2,3,4$, the direction indicator for the arc $L_i$ pointing to  $D_{i+1}$ 
from $D_i$ and there exist an inner edge of label $m$ 
for $D_i$ inward at $v_i$ and 
an inner edge of label $m$ for $D_{i+1}$ 
outward at $v_i$,
\end{enumerate}
here we understand $D_5=D_1$ 
(see Fig.~\ref{Fig09}(a)).
\end{enumerate}

Let $v$ be a white vertex of a chart and $e_1,e_2,e_3,e_4,e_5,e_6$ 
the edges situated clockwise around the vertex in this order. 
For each $i=1,2,3,4,5,6$, the edges $e_{i-1},e_{i+1}$ are called 
{\it side edges} of $e_i$, here we understand $e_0=e_6,e_7=e_1$. 
In Fig.~\ref{Fig09}(b), 
the edges $e_2$ and $e_6$ are side edges 
of $e_1$, also the edges $e_3$ and $e_5$ are the side edges of $e_4$.

\begin{lemma}
\label{PinWheelLemma}
{\rm (Pinwheel Lemma)}
If there exists a pinwheel, then we can decrease the number of crossings by $1$
and the number of white vertices by $4$.
Thus there does not exist a pinwheel in any minimal chart.
\end{lemma}

\begin{Proof}
Suppose that there exists a pinwheel in a minimal chart. 
For the pinwheel, we use the notation in 
Fig.~\ref{Fig09}(a).
For each $i=1,2,3,4$, 
let $v'_i$ be the closest vertex to the crossing on the arc $L_i$.
There exists an edge $e_i$ with vertices $c$ and $v'_i$. 
By direction indicators, 
the side edge of $e_i$ contained in $D_i$ (resp. $D_{i+1}$) is 
outward (resp. inward) at the vertex $v'_i$ 
(see Fig.~\ref{Fig09}(c)).
By applying four C-I-M2 moves we can get a mal-cycle 
(see Fig.~\ref{Fig09}(d)). 
This contradicts Lemma~\ref{MalCycleLemma}. 
\end{Proof}



\section{Consecutive Triplet Lemma}
\label{s:ConsecutiveTripletLemma}

In this section, we review a useful lemma
called Consecutive Triplet Lemma.

\begin{figure}
\begin{center}
\includegraphics{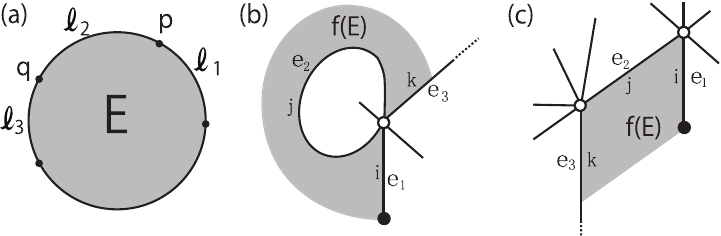}
\end{center}
\caption{ \label{Fig10} (b),(c) Consecutive triplets.}
\end{figure}

Let $E$ be a disk, and
$\ell_1,\ell_2,\ell_3$ three arcs on $\partial E$
such that each of $\ell_1\cap \ell_2$ and $\ell_2\cap \ell_3$ is one point and $\ell_1\cap \ell_3=\emptyset$
(see Fig.~\ref{Fig10}(a)),
say $p=\ell_1\cap \ell_2$,
$q=\ell_2\cap \ell_3$.
Let $\Gamma$ be a chart in a disk $D^2$.
Let $e_1$ be a terminal edge of 
 $\Gamma$. 
A triplet $(e_1,e_2,e_3)$ with distinct edges of $\Gamma$
is called 
a {\it consecutive triplet}
if there exists
a continuous map $f$ from the disk $E$ 
to the disk $D^2$ such that (see Fig.~\ref{Fig10}(b) and (c))
\begin{enumerate}
\item[(i)] the map $f$ is injective on $E-\{p,q\}$,
\item[(ii)] 
$f(\ell_3)$ is an arc in $e_3$, and $f({\rm Int}(E))\cap\Gamma=\emptyset$,
$f(\ell_1)=e_1$,
$f(\ell_2)=e_2$,
\item[(iii)]
each of $f(p)$ and $f(q)$ is a white vertex.
\end{enumerate}
If the label of $e_3$ is different
from that of $e_1$ 
then the consecutive triplet is said to be
{\it admissible}.


\begin{remark}
{\rm Let $(e_1,e_2,e_3)$
be a consecutive triplet. 
Since $e_2$ is an edge of $\Gamma$, 
the edge $e_2$ cannot contain a crossing.}
\end{remark}

\begin{lemma}
\label{ConsecutiveTripletLemma} 
{\rm (Consecutive Triplet Lemma)}
{\rm $($\cite[Lemma 1.1]{OneCrossing}, \cite[Lemma 3.1]{MinimalChart}$)$}
Any consecutive triplet 
in a minimal chart is admissible.
\end{lemma}

\begin{figure}
\begin{center}
\includegraphics{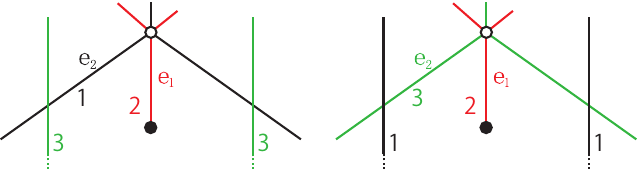}
\end{center}
\caption{ \label{TerminalEdgeLabel2} A terminal edge $e_1$ of label~$2$.}
\end{figure}

\begin{lemma}
\label{ConsecutiveTripletLemma4Chart}
Let $\Gamma$ be a minimal $4$-chart.
If a terminal edge $e_1$ is of label $2$,
then any side edge of $e_1$ contains a crossing
$($see Fig.~\ref{TerminalEdgeLabel2}$)$.
\end{lemma}

\begin{Proof}
Let $e_2$ be a side edge of $e_1$. 
Since $e_1$ is an edge of label~$2$,
the label of $e_2$ is $1$ or $3$.
Let $v$ be the vertex in $e_2$ but not contained in $e_1$.

If $v$ is a black vertex,
then $e_2$ is a terminal edge.
One of two terminal edges $e_1,e_2$ is not middle at the white vertex in $e_2$.
Thus we can apply a C-III move and
we eliminate the white vertex in $e_2$.
This contradicts the minimality of $\Gamma$.
 
If $v$ is a white vertex, 
then $v$ is contained in 
$\Gamma_1\cap \Gamma_2$ or $\Gamma_2\cap \Gamma_3$,
because $\Gamma$ is a $4$-chart.
Thus there exists an edge $e_3$ of label~$2$
such that the consecutive triplet $(e_1,e_2,e_3)$ is not admissible.
This contradicts Consecutive Triplet Lemma.
Therefore $v$ is a crossing.
\end{Proof}

The next lemma is a direct result of Consecutive Triplet Lemma.
This lemma will be used for terminal edges of label~$1$ or $3$ in minimal 4-charts.

\begin{figure}[bht]
\begin{center}
\includegraphics{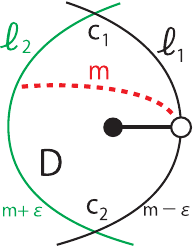}
\end{center}
\caption{ \label{Fig11} 
Any side edge of the terminal edge 
must intersect the arc $\ell_2$.}
\end{figure}

\begin{lemma}
\label{BridgeLemma}
Let $\Gamma$ be a minimal chart, and $m$ a label of the chart. 
Let $\varepsilon\in\{+1,-1\}$.
Let $c_1,c_2$ be crossings of the chart, and 
$D$ a disk whose boundary consists of two arcs $\ell_1,\ell_2$ 
of label $m-\varepsilon, m+\varepsilon$ respectively
with $\ell_1\cap\ell_2$ two crossings $c_1,c_2$.
If the disk $D$ does not contain any crossing except $c_1,c_2$, 
and if 
the interior of the disk $D$ does not contain any white vertex, 
then for any terminal edge of label $m-\varepsilon$ in the disk $D$, 
any side edge of the terminal edge
must intersect the arc $\ell_2$ 
$($see Fig.~\ref{Fig11}$)$.
\end{lemma}

\section{Anacanthous bodies $AB(\Gamma_1),AB(\Gamma_3)$}
\label{s:AnacanthousBody}

In this section, we investigate 
the set $\Gamma_1\cup \Gamma_3$ for a linear minimal 4-chart $\Gamma$.

Let $\Gamma$ be a linear minimal 4-chart.
Let $T(\Gamma)$ be the union of all the terminal edges of the chart $\Gamma$, 
$AB(\Gamma_1), AB(\Gamma_3)$ the closures of the sets 
$\Gamma_1-T(\Gamma), \Gamma_3-T(\Gamma)$ respectively.
We call $AB(\Gamma_1), AB(\Gamma_3)$ 
{\it anacanthous bodies} of label $1$ and label $3$ of $\Gamma$, respectively.

Let $\Gamma$ be a 4-chart.
For each label $k=1,3$, 
the closure of a connected component of $AB(\Gamma_k)-{\rm Cross}(\Gamma)$ is called 
an {\it  $AB$-component} of label $k$.

\begin{lemma} \LABEL{NoTree}
Let $\Gamma$ be a linear minimal $4$-chart without hoops nor free edges. 
Then for each $i=1, 3$, 
any connected component of $\Gamma_i$ contains 
a crossing. 
\end{lemma}

\begin{Proof}
Suppose that there exists a connected component $X$ 
of $\Gamma_i$ not containing a crossing. 
Then $X$ is also a connected component of $\Gamma_i-{\rm Cross}(\Gamma)$.
Since $\Gamma$ is linear,
the component $X$ is a tree. 
Since $\Gamma$ does not contain any free edges,
 $X$ is not a free edge.
Hence $X$ contains a white vertex. 
Let $Y$ be the tree obtained by omitting all the terminal edges 
from $X$. Let $v$ be an endpoint of $Y$. 
Then there exist two terminal edges at $v$. 
One of them is not middle at $v$. 
Applying a C-III move we can eliminate the vertex $v$. 
This contradicts the minimality of the chart.
\end{Proof}

\begin{lemma}
\label{DTLemma}
Let $\Gamma$ be a linear minimal $4$-chart with at most four crossings
and without hoops nor free edges. 
For $k=1,3$, 
\begin{enumerate}
\item[{\rm (a)}] any cycle in an $AB$-component of label $k$ 
contains at least two crossings,
\item[{\rm (b)}] any $AB$-component of  label $k$ 
does not contain any cycle. 
\end{enumerate}
\end{lemma}

\begin{figure}[h]
\begin{center}
\includegraphics{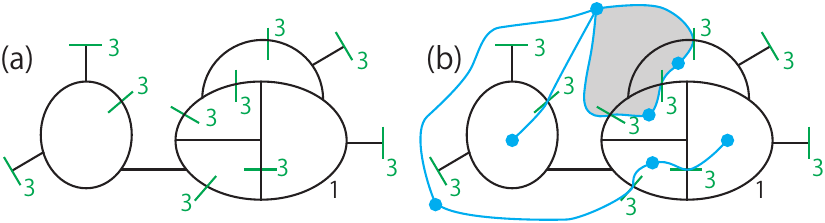}
\end{center}
\caption{\label{Fig12} 
(a) An $AB$-component $L$ of label~$1$.
(b) The graph $G_L$.
Here green short arcs are arcs in edges of label~$3$.}
\end{figure}

\begin{Proof}
We shall show Statement~(a). 
We show only for the case $k=1$. 
Let $L$ be an $AB$-component of label $1$. 
Let $C$ be a cycle in $L$.
Since $\Gamma$ is linear,
$AB(\Gamma_1)-{\rm Cross}(\Gamma)$ does not contain any cycle.
Thus $C$ contains at least one crossing.

Suppose that $C$ contains only one crossing $c_0$. 
Let $D$ be a disk  with $\partial D=C$. 

We claim that ${\rm Int}(D)$ contains at least two crossings.   
Let $L^*=D\cap AB(\Gamma_3)$. 
Then $L^*$ is not a terminal edge. 
Thus $L^*$ contains a white vertex. 
If $L^*$ is a tree,  
let $w$ be an endpoint of the tree different from the crossing. 
Then there exist two terminal edges at $w$ one of which is not middle at $w$. 
Hence we can eliminate the white vertex $w$ by a C-III move. 
This contradicts the minimality of the chart.

Hence $L^*$ contains a cycle $C^*$ contained in ${\rm Int}(D)$. 
Since the chart is linear, the cycle $C^*$ contains a crossing. 
If $C^*$ contains two crossings, then our claim holds. 
Suppose that the cycle contains only one crossing.
Let $D^*$ be the disk in ${\rm Int}(D)$ bounded by $C^*$. 
Let $L^{**}=D^*\cap AB(\Gamma_1)$. 
Then by the same way as in the proof of $L^*$,
the graph $L^{**}$ contains at least one crossing in ${\rm Int}(D^*)$. 
Hence ${\rm Int}(D)$ contains at least two crossings. 
Thus our claim holds.

Now by the same way as in the proof of $D$, 
the exterior of $D$ contains at least two crossings. 
Hence the chart $\Gamma$ contains at least five crossings. 
This contradicts the assumption for the number of crossings. 
Hence (a) holds. 

We shall show Statement~(b). 
Suppose that 
an $AB$-component $L$ of the chart $\Gamma$ contains a cycle
(see Fig.~\ref{Fig12}(a)). 
Then there are at least two complementary domains. 
Let 
$\mathcal{S}(L)$ be the set of all crossings each of which is on a cycle in $L$. 
For each complementary domain $D$ of $L$ in the sphere, 
we assign a point $p(D)$ in $D$. 
For each crossing $c$ in $\mathcal{S}(L)$, 
there are exactly two complementary domains $D_1,D_2$ 
each of whose closures contains the crossing $c$. 
Let $\ell(c)$ be a simple arc connecting $p(D_1)$ and $p(D_2)$ 
passing through the crossing $c$. 
We can assume that
for two different crossings $c_1,c_2$ in $\mathcal{S}(L)$, 
the simple arcs $\ell(c_1)$ and $\ell(c_2)$ 
do not intersect each other at their interiors.
Let $G_L=\cup\{\ell(c)~|~c\in \mathcal{S}(L)\}$, and $\mathcal{D}(L)$ the set of all complementary domains of $L$ in the sphere. 
Then we can consider $G_L$ as a graph whose edges are $\{\ell(c)~|~c\in \mathcal{S}(L)\}$ and whose vertices are $\{p(D) ~|~ D\in \mathcal{D}(L),$~${\rm Cl}(D)-D~{\rm contains~a~crossing~ in~}\mathcal{S}(L)\}$ 
(see Fig.~\ref{Fig12}(b)). 
Since $\mathcal{S}(L)\neq\emptyset$, we have $G_L\neq\emptyset$. 
If $G_L$ contains a cycle, then $L$ is not connected
(see the gray region in Fig.~\ref{Fig12}(b)).  
This is a contradiction. 
Thus $G_L$ is a tree. 
Hence the graph $G_L$ contains two endpoints $p(D'),p(D^*)$.
We can assume $D'$ does not contain the point at infinity of the sphere. 
If $p(D')$ is an  endpoint of $\ell(c)$, 
then the boundary ${\rm Cl}(D')-D'$ contains a cycle with only one crossing $c$.
This contradicts (a).
\end{Proof}

Let $\Gamma$ be a chart. 
Suppose that an object consists of 
some edges of $\Gamma$, arcs in edges of 
$\Gamma$ and arcs around white vertices.
Then the object is called a {\it pseudo-chart}.

\begin{figure}[h]
\begin{center}
\includegraphics{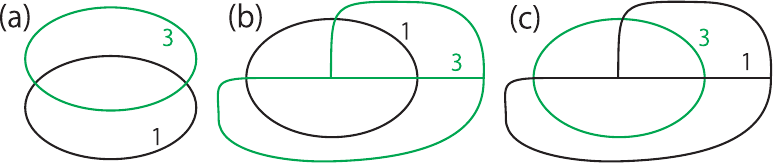}
\end{center}
\caption{\label{Fig13} 
 $AB(\Gamma_1) \cup AB(\Gamma_3)$.}
\end{figure}

\begin{lemma}
\label{Start} 
Let $\Gamma$ be a linear minimal $4$-chart with at most three crossings 
and without hoops nor free edges.
Then as the set, 
$AB(\Gamma_1) \cup AB(\Gamma_3)$ is homeomorphic to one of the three pseudo-charts drawn in Fig.~\ref{Fig13}. 
\end{lemma}

\begin{figure}[h]
\begin{center}
\includegraphics{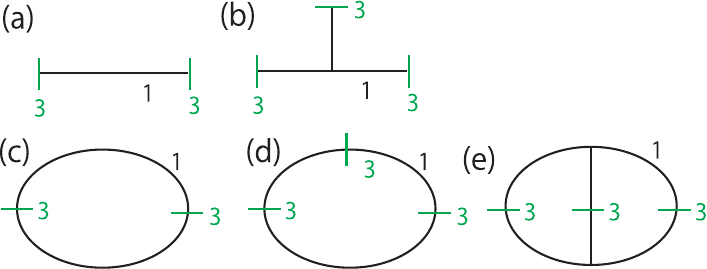}
\end{center}
\caption{\label{Fig14} 
(a),(b) $AB$-components of label $1$.
(c),(d),(e) $AB(\Gamma_1)$.
Here green short arcs are arcs in edges of label~$3$.}
\end{figure}

\begin{Proof}
Let $G$ be an $AB$-component of label $1$.
First, we shall show that $G$ is homeomorphic to one of the two pseudo-charts as drawn in Fig.~\ref{Fig14}(a),(b).

By Lemma~\ref{NoTree},
the $AB$-component $G$ contains at least one crossing.
By Lemma~\ref{DTLemma}(b),
 each crossing in $G$ is a vertex of degree 1 as a vertex of $G$. 
Moreover, each white vertex $w$  in $G$ is a vertex of degree 2 or 3
as a vertex of $G$.
If not, 
then there exist two terminal edges of label $1$ at $w$ one of which is not middle at $w$. 
Hence we can eliminate the white vertex by a C-III move. 
This contradicts the minimality of the chart.
Thus each white vertex in $G$ is a vertex of degree 2 or 3
as a vertex of $G$.

Since $\Gamma$ contains at most three crossings, 
the set $G$ contains at most three crossings.
Hence by Lemma~\ref{DTLemma}(b),
the set $G$ is a tree with at most three vertices of degree~$1$.
Thus $G$ is one of the two pseudo-charts as drawn in Fig.~\ref{Fig14}(a),(b).

If $G$ is the pseudo-chart as shown in Fig.~\ref{Fig14}(a),
then $AB(\Gamma_1)$ is a simple closed curve with at most three crossings
 as shown in Fig.~\ref{Fig14}(c),(d).
If $G$ is the pseudo-chart as shown in Fig.~\ref{Fig14}(b),
then $AB(\Gamma_1)$ is a $\theta$-curve as shown in
Fig.~\ref{Fig14}(e).

Similarly, the set $AB(\Gamma_3)$ is homeomorphic to one of the three pseudo-charts as shown in Fig.~\ref{Fig14}(c),(d),(e) by exchanging labels.
Therefore, $AB(\Gamma_1)\cup AB(\Gamma_3)$ is homeomorphic to one of the three pseudo-charts drawn in Fig.~\ref{Fig13}. 
\end{Proof}

\section{IO-paths}
\label{s:IO-Path}

In this section, we investigate simple arcs in  $(AB(\Gamma_1)\cup AB(\Gamma_3))-{\rm Cross}(\Gamma)$.

\begin{fact}
\label{Fact1}
{\rm Let $e$ be an edge of a linear minimal 4-chart $\Gamma$ not contained in 
$AB(\Gamma_1)\cup AB(\Gamma_3)$. 
If the edge $e$ is of label $1$ or $3$,
then the edge $e$ is a terminal edge.}
\end{fact}

\begin{figure}[htb]
\begin{center}
\includegraphics{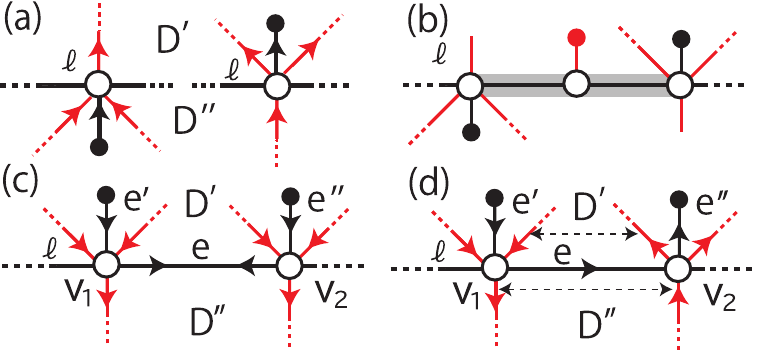}
\end{center}
\caption{\label{Fig15} 
(a) Fact~\ref{Fact2}.
(b) There exists a non-admissible consecutive triplet.
(c) Terminal edges $e',e''$ are contained in $D'$.
(d) The double-headed dotted arrows indicate the sites where C-I-M2 moves will be applied.}
\end{figure}

\begin{fact}
\label{Fact2}
{\rm Let $\Gamma$ be a linear minimal 4-chart.
Let $\ell$ be an anacanthous path of label~$1$ or $3$.
 Let $D'$ and $D''$ be disks with $D' \cap D''=\ell$.
If an inner edge for $D'$ is oriented outward $($resp. inward$)$ 
at an inner vertex $v$ of $\ell$, then
\begin{enumerate} 
\item[{\rm (1)}] any inner edge for the disk $D'$ containing the vertex $v$ 
is oriented outward $($resp. inward$)$ at $v$, and
\item[{\rm (2)}] any inner edge for the disk $D''$ containing the vertex $v$ 
is oriented inward $($resp. outward$)$ at $v$ $($see Fig.~\ref{Fig15}$($a$))$.
\end{enumerate} }
\end{fact}

Fact~\ref{Fact2} follows from Fact~\ref{Fact1}, Assumption~\ref{AssumeTerminal} and the definition of chart around a white vertex.

\begin{fact}
\label{Fact4}
{\rm Let $\Gamma$ be a linear minimal 4-chart.
Let $\ell$ be an anacanthous path of label $1$ or $3$.
If the arc $\ell$ contains at least two white vertices, then any edge of label~$2$ intersecting the interior of $\ell$ is not a terminal edge.
In other words, 
if a terminal edge of label~$2$ intersects the interior of the arc $\ell$, 
then the arc $\ell$ contains exactly one white vertex,
and the endpoints of $\ell$ are crossings (see Fig.~\ref{TerminalEdgeLabel2}).}
\end{fact}

Fact~\ref{Fact4} follows from Consecutive Triplet Lemma and Lemma~\ref{ConsecutiveTripletLemma4Chart} $($see Fig.~\ref{Fig15}$($b$))$.

\begin{fact}
\label{Fact5}
{\rm 
Let $\Gamma$ be a linear minimal 4-chart.
Let $\ell$ be an anacanthous path of label~$1$ or $3$.
 Let $D'$ and $D''$ be disks with $D' \cap D''=\ell$.
Let $e$ be an edge contained in the interior of the arc $\ell$, and 
$e', e''$ distinct terminal edges of label $1$ or $3$ 
each of which contains a vertex of the edge $e$. 
Then we have the following:
\begin{enumerate} 
\item[{\rm (1)}]
If $e'\subset D'$ then $e''\subset D''$. 
\item[{\rm (2)}]
If $e'\subset D''$ then $e''\subset D'$.
\end{enumerate} }
\end{fact}

\begin{Proof}
Suppose that $e', e''\subset D'$. 
If each of the edges $e',e''$ is oriented inward (resp. outward) at the vertex of the edge $e$, 
then the edge $e$ is oriented outward 
(resp. inward) at both white vertices of the edge $e$ 
(see Fig.~\ref{Fig15}(c)). 
This is a contradiction. 
If one of the edges $e',e''$ is oriented inward at a vertex of the edge $e$ and the other outward at the other vertex of the edge,
then we can eliminate the two white vertices of the edge $e$ 
by performing two C-I-M2 moves and one C-I-M3 move 
(see Fig.~\ref{Fig15}(d)).
\end{Proof}

Let $\ell$ be an arc consisting of edges of label $m$ in a chart $\Gamma$ and  
$D',D''$ disks with $D'\cap D''=\ell$.
If, for each inner vertex $v$ of the arc $\ell$, 
any inner edge for $D'$ containing the vertex $v$ 
is outward at the vertex $v$, and 
any inner edge for $D''$ containing the vertex $v$ 
is inward at the vertex $v$,
then we call the arc $\ell$ an {\it IO-path} of label $m$, 
and 
we put a direction indicator on $\ell$ pointing to $D'$ from $D''$.

\begin{note}
{\rm We regard an arc with no inner vertices as an IO-path.
In this case, we can put a direction indicator for the arc 
as we wish.}
\end{note}

\begin{lemma}
\label{IO-pathLemma}
Let $\Gamma$ be a linear minimal $4$-chart.
Let $\ell$ be an anacanthous path of label~$1$ or $3$.
Then the arc $\ell$ is an IO-path.
\end{lemma}

\begin{figure}[htb]
\begin{center}
\includegraphics{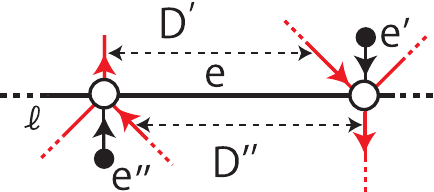}
\end{center}
\caption{\label{Fig16} 
The double-headed dotted arrows indicate the sites where
C-I-M2 moves will be applied.}
\end{figure}

\begin{Proof}
Let $D',D''$ be disks with $D'\cap D''=\ell$. 
If $\ell$ is not an IO-path, then, by the help of Fact~\ref{Fact5}, 
there exist an edge $e$ and two terminal edges $e',e''$ such that
\begin{enumerate}
\item[(1)] 
the edge $e$ is contained in the interior of the arc $\ell$,
\item[(2)] 
the terminal edge $e'$ is an inner edge for $D'$ inward 
(resp. outward) at a white vertex of $e$, and
the terminal edge $e'' $ is an inner edge for $D''$ inward
(resp. outward) at the other white vertex of $e$.
\end{enumerate}
But we can eliminate the two white vertices of $e$ 
by two C-I-M2 moves and one C-I-M3 move 
(see Fig.~\ref{Fig16}).
This contradicts the minimality of the chart.
\end{Proof}

\begin{figure}[htb]
\begin{center}
\includegraphics{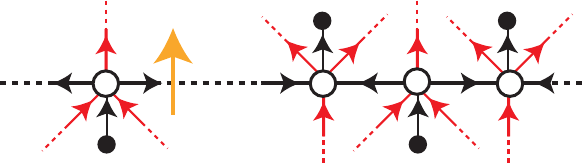}
\end{center}
\caption{\label{FigAnacanthousPath} 
An anacanthous path is an IO-path.}
\end{figure}

By Fact~\ref{Fact5} and Lemma~\ref{IO-pathLemma},
in a linear minimal $4$-chart,
an anacanthous path of label 1 or 3
is an IO-path as shown in Fig.~\ref{FigAnacanthousPath}.

\section{Linear minimal 4-charts with two crossings}
\label{s:4ChartTwoCrossing}

In this section, we shall show that 
after deleting free edges and hoops, any linear minimal 4-chart with two crossings
is lor-equivalent to the 4-chart describing a turned torus-link of a Hopf link.

\begin{figure}[htb]
\begin{center}
\includegraphics{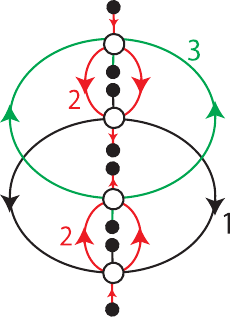}
\end{center}
\caption{\label{Fig17} 
The linear minimal 4-chart with two crossings.}
\end{figure}

\begin{proposition} 
\label{MainTheorem2} 
After deleting free edges and hoops,
any linear minimal $4$-chart with two crossings is 
lor-equivalent to the chart as shown in Fig.~\ref{Fig17}. 
\end{proposition}

\begin{Proof}
Let $\Gamma$ be a linear minimal $4$-chart with two crossings.
By Lemma~\ref{4ChartSimpleHoop},
any hoop is simple.
Thus by Assumption~\ref{NoSimpleHoop},
all free edges and hoops are moved into a small neighbourhood of the point at infinity
$\infty$.
Therefore, we can assume that $\Gamma$ contains neither free edges nor hoops.

By Lemma~\ref{Start},
the anacanthous body $AB(\Gamma_1)\cup AB(\Gamma_3)$ consists of four $AB$-components $\ell_1,\ell_2,\ell_3,\ell_4$ 
where $\ell_1,\ell_3$ are arcs of label~$1$
and $\ell_2,\ell_4$ are arcs of label~$3$ 
as shown in Fig.~\ref{Fig18}(a).
Let $D_1,D_2,D_3,D_4$ be disks each of which is the closure of a complementary domain of $AB(\Gamma_1) \cup AB(\Gamma_3)$ on the sphere as shown in Fig.~\ref{Fig18}(b).
Let $c_1,c_2$ be crossings 
as shown in Fig.~\ref{Fig18}(c).

\begin{figure}[h]
\begin{center}
\includegraphics{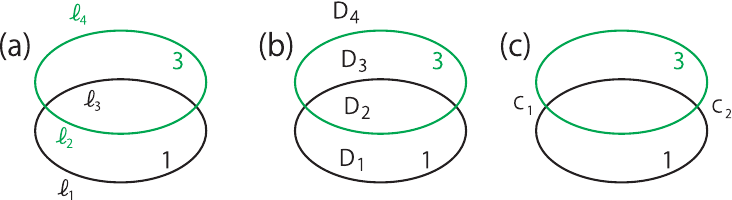}
\end{center}
\caption{\label{Fig18} 
(a) Arcs.
(b) Domains.
(c) Crossings. }
\end{figure}

Since all of $\ell_1,\ell_2,\ell_3,\ell_4$ 
are anacanthous paths,
we can put a direction indicator for each of the arcs 
$\ell_1,\ell_2,\ell_3,\ell_4$
by Lemma~\ref{IO-pathLemma}.
If necessary we reverse the orientation of all the edges,
we can assume that a direction indicator for the arc $\ell_1$
points to the disk $D_1$ from $D_4$
(see Fig.~\ref{Fig19}(a)).

{\bf Claim.} Each of the arcs $\ell_1,\ell_2,\ell_3,\ell_4$ contains at most  one inner vertex. 

{\it Proof of Claim.}
Suppose that the arc $\ell_1$ contains at least two inner vertices. Then by Fact~\ref{Fact5},
there exist at least three inner edges of label~$2$ for the disk $D_1$ each of which is oriented outward at an inner vertex of $\ell_1$
(see Fig.~\ref{Fig19}(b)).
Thus by Fact~\ref{Fact4},
none of the three edges are terminal edges.
Hence by considering the direction indicator for the arc $\ell_1$,
all of the three edges contain inner vertices of $\ell_2$.
Thus the arc $\ell_2$ contains at least two inner vertices
and the direction indicator for $\ell_2$ points to $D_2$ from $D_1$ (see Fig.~\ref{Fig19}(c)). 

Similarly, there exist at least three inner edges for $D_2$
each of which is oriented outward at an inner vertex of $\ell_2$ 
and oriented inward at an inner vertex of $\ell_3$.
Moreover, there exist at least three inner edges for $D_3$
each of which is oriented outward at an inner vertex of $\ell_3$ 
and oriented inward at an inner vertex of $\ell_4$.
Further, there exist at least three inner edges for $D_4$
each of which is oriented outward at an inner vertex of $\ell_4$ 
and oriented inward at an inner vertex of $\ell_1$.
Therefore, we can find a pinwheel around the crossing $c_1$
(see  Fig.~\ref{Fig19}(d)). 
This contradicts Lemma~\ref{PinWheelLemma}.
Hence the arc $\ell_1$ contains at most one inner vertex.

Similarly, we can show that each of $\ell_2,\ell_3,\ell_4$ 
contains at most one inner vertex.
Thus Claim holds.
\hfill{$\square$}

\begin{figure}[h]
\begin{center}
\includegraphics{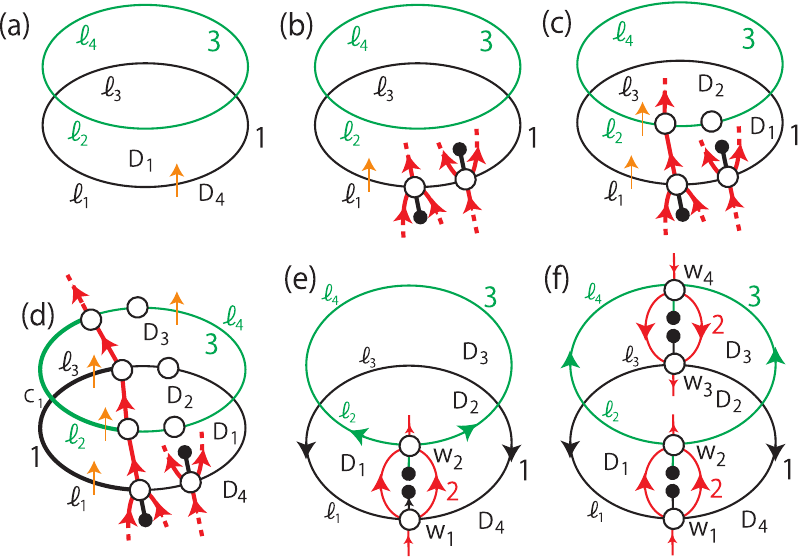}
\end{center}
\caption{\label{Fig19} 
(a) The direction indicator for $\ell_1$.
(b) The arc $\ell_1$ contains at least two inner vertices.
(c) The arc $\ell_2$ contains at least two inner vertices.
(d) A pinwheel around the crossing $c_1$.
(e) The arc $\ell_1$ contains one inner vertex.
(f) The arc $\ell_3$ contains one inner vertex.}
\end{figure}

Now we start the proof of Proposition~\ref{MainTheorem2}.

If none of $\ell_1,\ell_2,\ell_3,\ell_4$
contain inner vertices,
then we can eliminate the two crossings by a C-I-R2 move.
This contradicts the minimality of the chart.
Thus one of $\ell_1,\ell_2,\ell_3,\ell_4$ contains
at least one inner vertex.

Without loss of generality,
we can assume that
the arc $\ell_1$ contains an inner vertex $w_1$.
Let $e_1$ be the terminal edge of label $1$ containing $w_1$.
If necessary we move the point at infinity $\infty$ into $D_1$,
we can assume that the edge $e_1$ is an inner edge for the disk $D_1$.
Moreover, we can assume that the edge $e_1$ is oriented outward at $w_1$.

There are two side edges of label $2$ of $e_1$.
By Lemma~\ref{BridgeLemma},
both of the edges contain inner vertices of $\ell_2$.
Thus by Claim, both of the edges contain the same inner vertex
$w_2$ of $\ell_2$.
Hence we have the pseudo-chart as shown in Fig.~\ref{Fig19}(e).

If neither $\ell_3$ nor $\ell_4$ contains
an inner vertex,
then we can eliminate the two crossings by a C-I-R2 move.
This contradicts the minimality of the chart.
Thus one of $\ell_3,\ell_4$ contains
at least one inner vertex.

If necessary we change labels and we move the point at infinity $\infty$ into $D_2$,
we can assume that $\ell_3$ contains an inner vertex $w_3$.
Let $e_3$ be the terminal edge of label $1$ containing $w_3$.
Then $e_3\subset D_2$ or $e_3\subset D_3$.

If $e_3\subset D_2$,
then the side edges of label $2$ of $e_3$ 
contain inner vertices of $\ell_1$
by Lemma~\ref{BridgeLemma}.
Thus the arc $\ell_1$ must contain at least two inner vertices.
This contradicts Claim.
Hence  $e_3\subset D_3$.

Thus by Lemma~\ref{BridgeLemma} and Claim,
the side edges of $e_3$ contain the same inner vertex $w_4$
of $\ell_4$.
Therefore, we have
the pseudo-chart as shown in Fig.~\ref{Fig19}(f).
Moreover, we can see that
the two inner edges for $D_2$ are oriented outward at $w_2,w_3$,
respectively, and
the two inner edges for $D_4$ are oriented inward at $w_1,w_4$,
respectively.
Hence the four edges must be terminal edges.
Therefore, we have the chart as shown in Fig.~\ref{Fig17}.
\end{Proof}


\section{Proof of the main theorem}
\label{s:ProofTheorem}

In this section, we shall show the main theorem.

From now on, by Lemma~\ref{4ChartSimpleHoop} and Assumption~\ref{NoSimpleHoop},
we assume that
\begin{enumerate}
\item[] {\it the chart $\Gamma$ is a linear minimal $4$-chart with three crossings without free edges nor hoops.} 
\end{enumerate}

\begin{figure}[h]
\begin{center}
\includegraphics{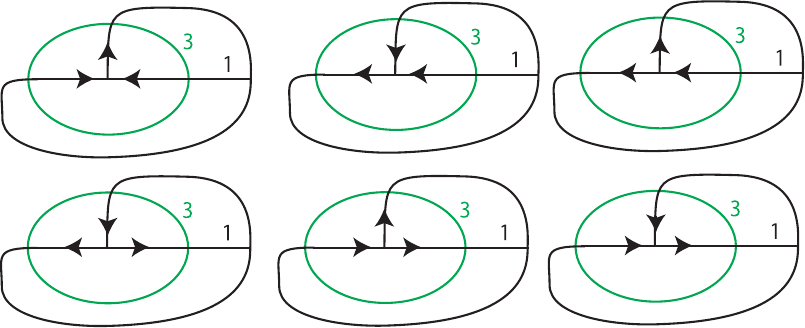}
\end{center}
\caption{\label{Fig20} 
The set $AB(\Gamma_1) \cup AB(\Gamma_3)$.}
\end{figure}

By exchanging labels, 
two pseudo-charts in Fig.~\ref{Fig13}(b),(c)
 are lor-equivalent,
we adopt Fig.~\ref{Fig13}(c) 
for $AB(\Gamma_1)\cup AB(\Gamma_3)$.
Further, since six pseudo-charts with orientation in 
Fig.~\ref{Fig20} 
are lor-equivalent, 
we adopt the first one in Fig.~\ref{Fig20}.

\begin{figure}[hbt]
\begin{center}
\includegraphics{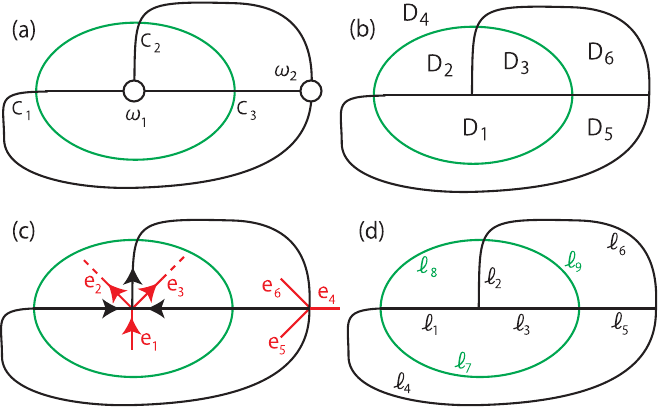}
\end{center}
\caption{\label{Fig21} 
(a) Vertices and crossings.
(b) Domains.
(c) Edges.
(d) Arcs.}
\end{figure}

Let $w_1,w_2$ be white vertices on $ AB(\Gamma_1)$, 
and $c_1,c_2,c_3$ crossings 
as shown in Fig.~\ref{Fig21}(a).
Let $D_1,D_2,D_3,D_4,D_5,D_6$ be disks each of which is the closure of a complementary domain of $AB(\Gamma_1) \cup AB(\Gamma_3)$ on the sphere 
as shown in Fig.~\ref{Fig21}(b).
Let $e_1,e_2,e_3$ be edges of label $2$ containing the white vertex $w_1$, and $e_4,e_5,e_6$ be edges of label $2$ containing the white vertex $w_2$ as shown in Fig.~\ref{Fig21}(c).
By orientation of $\Gamma_1$ around the vertex $w_1$, 
the edge $e_1$ is oriented inward at the vertex $w_1$ and 
$e_2,e_3$ oriented outward at the vertex $w_1$.
Now the set $AB(\Gamma_1)-\{w_1,w_2,c_1,c_2,c_3\}$ consists of 
six open arcs.
Let $\ell_1,\ell_2,\ell_3,\ell_4,\ell_5,\ell_6$ 
be the closures of the six open arcs 
where $\ell_1,\ell_2,\ell_3$ contain the vertex $w_1$, and 
$\ell_4,\ell_5,\ell_6$ the vertex $w_2$ 
as shown in Fig.~\ref{Fig21}(d).
Let $\ell_7,\ell_8,\ell_9$ be the closures of three open arcs 
$AB(\Gamma_3)-\{c_1,c_2,c_3\}$ as shown in Fig.~\ref{Fig21}(d).

\begin{fact}
\label{Fact3}
{\rm Each of the disks $D_1,D_2,D_3,D_4,D_5,D_6$ 
satisfies the two conditions for the disk in Lemma~\ref{BridgeLemma}, hence the result of Lemma~\ref{BridgeLemma}.}
\end{fact}

Since all of $\ell_1,\ell_2,\ell_3,\ell_4,\ell_5,\ell_6,\ell_7,\ell_8,\ell_9$ are anacanthous paths,
by Lemma~\ref{IO-pathLemma}
we can put a direction indicator for each of the arcs  
$\ell_1, \ell_2, \ell_3,\ell_4,\ell_5,\ell_6,\ell_7,\ell_8,\ell_9$.

\begin{figure}[h]
\begin{center}
\includegraphics{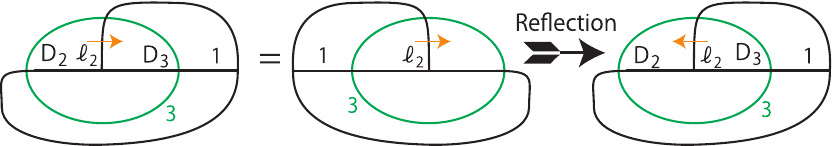}
\end{center}
\caption{\label{Fig22} 
Direction indicators for the arc $\ell_2$.}
\end{figure}

As shown in Fig.~\ref{Fig22}, 
the pseudo-chart with the arc $\ell_2$ 
having a direction indicator pointing to $D_2$ from $D_3$ 
is lor-equivalent to 
the pseudo-chart with the arc $\ell_2$ 
having a direction indicator pointing to $D_3$ from $D_2$.
Hence, we assume that the arc $\ell_2$ has 
a direction indicator pointing to $D_2$ from $D_3$.

The following will be shown in the next subsection.

\begin{proposition}
\label{E1Terminal}
The edge $e_1$ is a terminal edge.
\end{proposition}

\subsection{Case that the edge $e_1$ is not a terminal edge}

We shall prove Proposition~\ref{E1Terminal} by five claims.
Suppose the edge $e_1$ is not a terminal edge.

\begin{figure}[h]
\begin{center}
\includegraphics{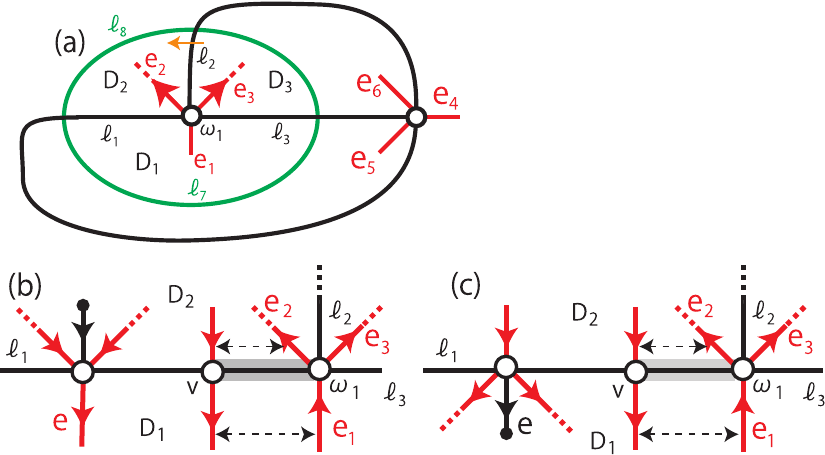}
\end{center}
\caption{\label{Fig23} 
(a) Neither $e_2$ nor $e_3$ are terminal edges.
(b),(c) Proof of Claim~\ref{Claim1}.
The double-headed dotted arrows indicate the sites where
C-I-M2 moves will be applied.}
\end{figure}

\begin{claim}
\label{Claim1} 
The direction indicator for the arc $\ell_1$
points to $D_2$ from $D_1$, and 
the direction indicator for the arc $\ell_3$
points to $D_3$ from $D_1$.
\end{claim}

{\it Proof of Claim~\ref{Claim1}.} 
We start from Fig.~\ref{Fig23}(a). 
Suppose that there exists an inner edge $e$ for $D_1$ outward 
at a vertex on the arc $\ell_1$. 
Let $v$ be the closest vertex to $w_1$ on the arc $\ell_1$. 
By Lemma~\ref{IO-pathLemma}, the arc $\ell_1$ is an IO-path. 
Hence any inner edge for the disk $D_1$ 
containing the vertex $v$ is 
oriented outward at $v$. 
Then by Fact~\ref{Fact2} there exist 
an inner edge for $D_1$ of label $2$ outward at $v$ 
and 
an inner edge for $D_2$ of label $2$ inward at $v$ 
(see Fig.~\ref{Fig23}(b),(c)).
Hence by the help of the edges $e_1$ and $e_2$, 
applying two C-I-M2 moves along the double-headed dotted lines 
and one C-I-M3 move we can eliminate two white vertices. 
This contradicts the minimality of the chart.
Thus, the direction indicator for the arc $\ell_1$ points 
to $D_2$ from $D_1$.

Similarly, 
we can show that the direction indicator for the arc $\ell_3$ 
points to $D_3$ from $D_1$.
Hence, Claim~\ref{Claim1} holds. \hfill {$\square$}\vspace{1.5em} 

\begin{figure}[h]
\begin{center}
\includegraphics{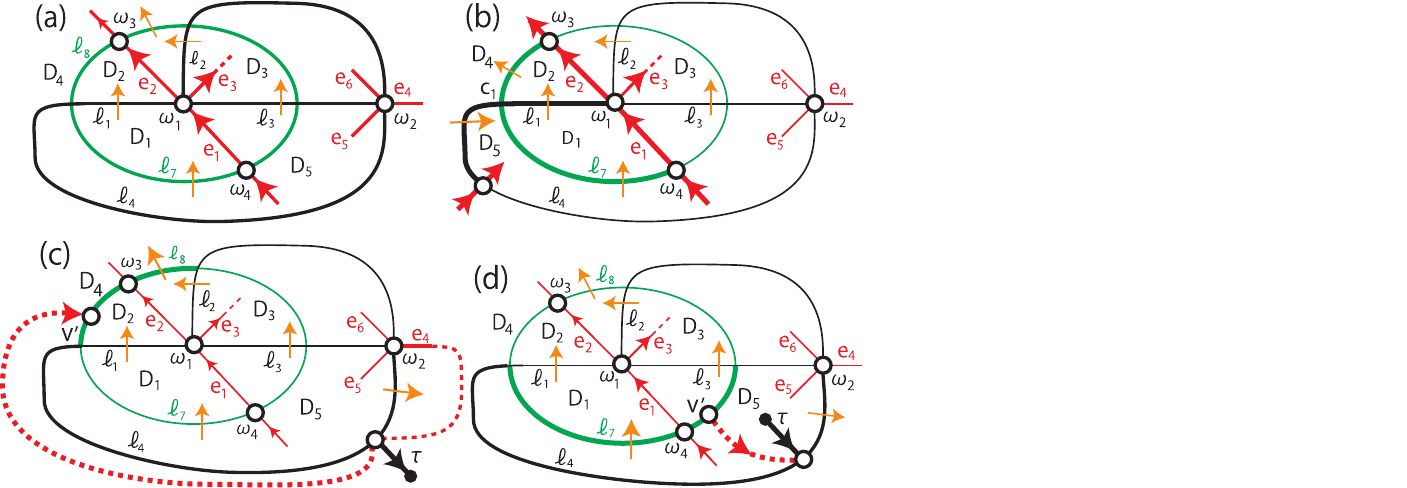}
\end{center}
\caption{\label{Fig24} 
(a) The five direction indicators for $\ell_1,\ell_2,\ell_3,\ell_7,\ell_8$.
(b) Case~$1$ of Claim~\ref{Claim2}. A pinwheel around the crossing $c_1$. 
(c) Case~$2.1$ of Claim~\ref{Claim2}. The terminal edge $\tau$ is an inner edge for $D_4$.
(d) Case~$2.2$ of Claim~\ref{Claim2}. The terminal edge $\tau$ is an inner edge for $D_5$.}
\end{figure}

Since the direction indicator for $\ell_2$ points to $D_2$ from $D_3$, 
Claim~$1$ implies that 
the edge $e_2$ must have a white vertex $w_3$ on the arc $\ell_8$. 
Thus, the direction indicator for the arc $\ell_8$ points to $D_4$ 
from $D_2$ (see Fig.~\ref{Fig24}(a)).

By Claim~\ref{Claim1}, 
we put direction indicators to the arcs $\ell_1,\ell_3$ 
as shown in Fig.~\ref{Fig24}(a). 
Hence, the edge $e_1$ contains a vertex $w_4$ on the arc $\ell_7$.
Thus, the direction indicator for the arc $\ell_7$ points to 
$D_1$ from $D_5$.

\begin{claim}
\label{Claim2} 
The arc $\ell_4$ does not contain any inner vertex.
\end{claim}

{\it Proof of Claim~\ref{Claim2}.}
We start from Fig.~\ref{Fig24}(a). 
Suppose that the arc $\ell_4$ contains an inner vertex $v$. 
There are two cases.

{\bf Case 1.} The direction indicator for the arc $\ell_4$ 
points to $D_5$ from $D_4$.

{\bf Case 2.} The direction indicator for the arc $\ell_4$ 
points to $D_4$ from $D_5$.

For Case~$1$. By Fact~\ref{Fact1},
there exist terminal edges of label~$3$ at white vertices $w_3,w_4$,
respectively.
Hence, we can find a pinwheel with the vertices $w_1,w_3,v,w_4$ around the crossing $c_1$ 
(see Fig.~\ref{Fig24}(b)). 
This contradicts Lemma~\ref{PinWheelLemma}.

For Case~$2$. By Fact~\ref{Fact1}, 
there exists a terminal edge $\tau$ of label~$1$ 
at the vertex $v$.  There are two cases.

{\bf Case 2.1.}
 The terminal edge $\tau$ is an inner edge for $D_4$.

{\bf Case 2.2.} The terminal edge $\tau$ is an inner edge for $D_5$.

For Case~$2.1$. 
One of the side edges of the terminal edge $\tau$ 
must contain a vertex $v'$ on the arc $\ell_8$ by Lemma~\ref{BridgeLemma}. 
The side edge must be inward at the vertex $v'$. But this is impossible because the direction indicator for the edge $\ell_8$ points to $D_4$ from $D_2$ (see Fig.~\ref{Fig24}(c)). 

For Case~$2.2$. 
One of the side edges of the terminal edge $\tau$ must contain a vertex $v'$ on the arc $\ell_7$ by Lemma~\ref{BridgeLemma}. The side edge must be outward at the vertex $v'$. But this is impossible because the direction indicator for the edge $\ell_7$ points to $D_1$ from $D_5$ 
(see Fig.~\ref{Fig24}(d)).
Hence, Claim~\ref{Claim2} holds. 
\hfill {$\square$}\vspace{1.5em}

\begin{figure}[h]
\begin{center}
\includegraphics{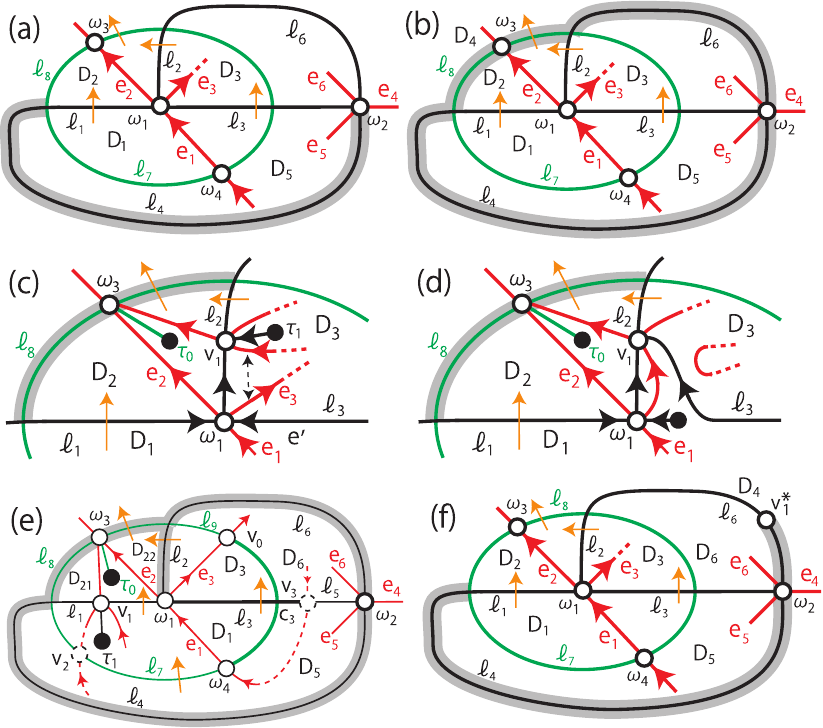}
\end{center}
\caption{\label{Fig25} 
(a) The arc $\ell_4$ does not contain any inner vertex.
(b) Proof of Claim~\ref{Claim3}. The arc $\ell_6$ does not contain any inner vertex,
and the arc $\ell_8$ does not contain any inner vertex except $w_3$. 
(c),(d) Case~2.1 of Claim~\ref{Claim3}. 
(e) Case~2.2 of Claim~\ref{Claim3}. 
(f) The arc $\ell_6$ contains an inner vertex.}
\end{figure}

Now we have the figure in Fig.~\ref{Fig25}(a).

\begin{claim}
\label{Claim3} 
There exists an inner vertex of the arc $\ell_6$.
\end{claim}

{\it Proof of Claim~\ref{Claim3}.}
We start from Fig.~\ref{Fig25}(a).
Suppose that there does not exist any inner vertex of the arc $\ell_6$. 
Then we have
\begin{enumerate}
\item[(1)]  the edge $e_4$ is only the edge of label $2$ 
containing a white vertex on the arc $\ell_4\cup \ell_6$.
\end{enumerate}
Hence, we have
\begin{enumerate}
\item[(2)] there does not exist an inner vertex of $\ell_8$ except $w_3$ (see Fig.~\ref{Fig25}(b)).
\end{enumerate}

By Fact~\ref{Fact1}, 
there exists a terminal edge $\tau_0$ of label $3$ containing the vertex $w_3$ 
on the arc $\ell_8$.
There are two cases. 

{\bf Case~1.} The terminal edge $\tau_0$ is an inner edge for the disk $D_4$.

{\bf Case~2.} The terminal edge $\tau_0$ is an inner edge for the disk $D_2$.

For Case~$1$. The two side edges of the terminal edge $\tau_0$  are edges
of label~$2$. The two side edges must have white vertices on the arc 
$\ell_4\cup \ell_6$ by Lemma~\ref{BridgeLemma}. 
Namely, there must exist at least two white vertices on the arc $\ell_4\cup \ell_6$. 
This contradicts~(1).

For Case~$2$. The edge $e_2$ splits the disk $D_2$ into two disks. 
Let $D_{21}$ be one of the two disks which contains the arc $\ell_1$, 
and $D_{22}$ the disk which contains the arc $\ell_2$. 
There are two cases.

{\bf Case~2.1.} The terminal edge $\tau_0$ is an inner edge for the disk $D_{22}$.

{\bf Case~2.2.} The terminal edge $\tau_0$ is an inner edge for the disk $D_{21}$.

For Case~$2.1$. The side edge of $\tau_0$ different from $e_2$ contains 
a white vertex $v_1$ on the arc $\ell_2$ by Lemma~\ref{BridgeLemma}. 
Let $\tau_1$ be the terminal edge of label~$1$ containing the vertex. 
By (2) and Lemma~\ref{BridgeLemma}, 
the terminal edge $\tau_1$ is an inner edge for $D_3$ and 
there does not exist a vertex between $v_1$ and $w_1$ on $\ell_2$. 
Since the direction indicator for the arc $\ell_2$ points to $D_2$ 
from $D_3$, the terminal edge $\tau_1$ is inward at the vertex $v_1$, 
and so are the side edges of the terminal edge $\tau_1$ 
(see Fig.~\ref{Fig25}(c)). 
Let $e'$ be the edge of label~$1$ in $\ell_3$ containing 
the vertex $w_1$. 
If necessary, by applying a C-I-M2 move 
between $e_3$ and a side edge of $\tau_1$, 
and then by applying a C-I-M2 move
between the terminal edge $\tau_1$ and the edge $e'$, 
the edge $e'$ changes to a terminal edge at $w_1$ 
(see Fig.~\ref{Fig25}(d)). 
Thus, we can eliminate the white vertex $w_1$ by a C-III move.
This contradicts the minimality of the chart.

For Case~$2.2$.  
First, we shall show that $\ell_2$ does not contain any inner vertex. If there exists an inner vertex $v$ of $\ell_2$,
then there exists an inner edge of label~$2$ 
for $D_{22}$ containing $v$.
By Fact~\ref{Fact4},  
the inner edge is not a terminal edge.
Hence considering the direction indicator for $\ell_2$,
the inner edge contains an inner vertex of $\ell_8$ different from $w_3$.
However, this contradicts (2).
Thus $\ell_2$ does not contain any inner vertex.

Since the edge $e_3$ is not a terminal edge,
the edge $e_3$ must contain an inner vertex $v_0$ for $\ell_9$
by considering the direction indicator for $\ell_3$
(see Fig.~\ref{Fig25}(e)).

Now, by Lemma~\ref{BridgeLemma}, 
the side edge of $\tau_0$ different from $e_2$ 
contains a white vertex $v_1$ on the arc $\ell_1$ 
at which the edge is outward. 
Let $\tau_1$ be the terminal edge of label $1$ 
containing the vertex $v_1$. 
Since there does not exist any inner vertex 
on the arc $\ell_8\cap D_{21}$ by (2), 
the terminal edge $\tau_1$ is 
an inner edge for the disk $D_1$ (see Fig.~\ref{Fig25}(e)). 
Again by Lemma~\ref{BridgeLemma}, 
the far side edge of the terminal edge $\tau_1$ from 
the edge $e_1$, contains 
a white vertex $v_2$ on the arc $\ell_7$ different from $w_4$.
Since there exist two white vertices on the arc $\ell_7$, 
each of the inner edges for the disk $D_5$ of label $2$ is not 
a terminal edge by Fact~\ref{Fact4}, and is 
inward at a vertex on the arc $\ell_7$ 
(see Fig.~\ref{Fig25}(e)). 
Thus, there exists a white vertex $v_3$ on $\ell_5$ 
at which one of the two inner edges for $D_5$ is outward.
Hence, there exists an inner edge of label~$2$ for $D_5$ outward at $v_3$ 
and an inner edge of label~$2$ for $D_6$ inward at $v_3$.
Thus, we can find a pinwheel with vertices $w_1,v_0,v_3,w_4$ 
around the crossing $c_3$
(see Fig.~\ref{Fig25}(e)).
This contradicts Lemma~\ref{PinWheelLemma}.

Therefore Claim~\ref{Claim3} holds.\hfill {$\square$}\vspace{1.5em} 


Let $v^*_1$ be the white vertex on the arc $\ell_6$, assured by Claim~$3$,
closest to the vertex $w_2$ on the arc $\ell_6$ 
(see Fig.~\ref{Fig25}(f)).
By Fact~\ref{Fact1}, 	
there exists a terminal edge $\tau$ of label $1$ containing 
the vertex $ v^*_1$.
Then, the terminal edge $\tau$ is an inner edge for the disk $D_4$, or 
an inner edge for the disk $D_6$.

\begin{claim}
\label{Claim4}
 The terminal edge $\tau$ is an inner edge of the disk $D_6$.
\end{claim}

{\it Proof of Claim~\ref{Claim4}.}
Suppose that the terminal edge $\tau$ is an inner edge for the disk $D_4$. 
There are two cases.

{\bf Case~1.} The terminal edge $\tau$ is outward at the vertex $v^*_1$.

{\bf Case~2.} The terminal edge $\tau $ is inward at the vertex $v^*_1$.

For Case~$1$. By Lemma~\ref{BridgeLemma}, each of the two side edges 
of $\tau $ of label~$2$ must contain a white vertex on the arc $\ell_8$. 
Since the terminal edge $\tau $ is outward at $ v^*_1$, each of 
the side edges must be inward at the white vertex on the arc $\ell_8$. 
But this is impossible because the direction indicator 
for the arc $\ell_8$ points to $D_4$ from $D_2$.

For Case~$2$. 
There are two cases.

{\bf Case~2.1.} The edge $e_4$ is inward at $w_2$.

{\bf Case~2.2.} The edge $e_4$ is outward at $w_2$.

Let $\tilde{e}$ be the edge contained in $\ell_6$ with two white vertices 
$w_2, v^*_1$. 
Since the terminal edge $\tau $ is inward at the vertex $v^*_1$,
 the edge $\tilde{e}$ is inward at $w_2$.
Thus, $e_5$ is outward at $w_2$.

\begin{figure}[h]
\begin{center}
\includegraphics{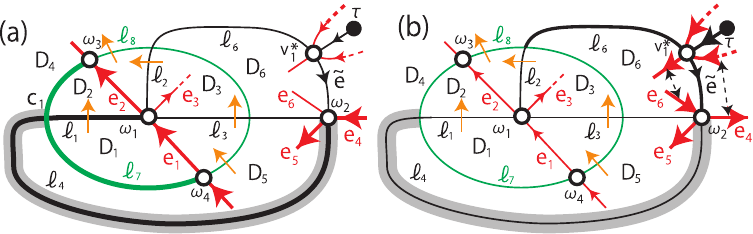}
\end{center}
\caption{\label{Fig26} 
The terminal edge $\tau$ is an inner edge for $D_4$.
(a) Case 2.1 of Claim~\ref{Claim4}. 
(b) Case 2.2 of Claim~\ref{Claim4}. }
\end{figure}

For Case~$2.1$. 
Since the edge $e_5$ is outward at $w_2$, 
we can find a pinwheel around the crossing $c_1$ containing 
the vertices $w_1,w_3,w_2,w_4$ (see Fig.~\ref{Fig26}(a)). 
This contradicts Lemma~\ref{PinWheelLemma}.

For Case~$2.2$. 
Since the edge $\tilde{e}$ is inward at $w_2$ and 
the edge $e_4$ is outward at $w_2$, 
the orientation of edges around the vertex $w_2$ is determined 
by the definition of the orientation around a white vertex
(see Fig.~\ref{Fig26}(b)). 
By performing two C-I-M2 moves and one C-I-M3 move,
we can eliminate the white vertices $w_2$ and $ v^*_1$. 
This contradicts the minimality of the chart.

Therefore, Claim~\ref{Claim4} holds.\hfill {$\square$}\vspace{1.5em}

Therefore, the terminal edge $\tau$ must be 
an inner edge for the disk $D_6$.

\begin{claim}
\label{Claim5}
 The terminal edge $\tau$ is outward at the vertex $ v^*_1$.
\end{claim}

{\it Proof of Claim~\ref{Claim5}.} 
Suppose that the terminal edge $\tau$ is inward at the vertex $ v^*_1$.
Then, the direction indicator for $\ell_6$ points to $D_4$ from $D_6$ 
by Lemma~\ref{IO-pathLemma}.
Let $\tilde{e}$ be the edge contained in $\ell_6$ with two white vertices 
$w_2, v^*_1$. 
Then the edge $\tilde{e}$ is inward at $w_2$ .
Thus, $e_5$ is outward at $w_2$.
There are two cases.

{\bf Case~1.} The edge $e_4$ is inward at the vertex $w_2$.

{\bf Case~2.} The edge $e_4$ is outward at the vertex $w_2$.
 
\begin{figure}[h]
\begin{center}
\includegraphics{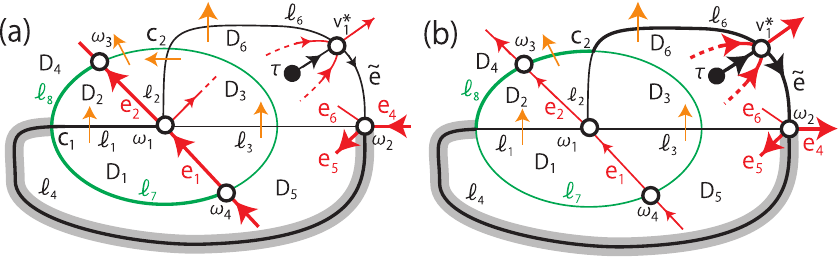}
\end{center}
\caption{\label{Fig27} 
The terminal edge $\tau$ is an inner edge for $D_6$,
and is inward at $v_1^*$.
(a) Case~$1$ of Claim~\ref{Claim5}. 
(b) Case~$2$ of Claim~\ref{Claim5}. }
\end{figure}

For Case~$1$. We have a pinwheel 
around the crossing $c_1$ containing 
the vertices $w_1,w_3,w_2,w_4$ as in the proof of Case~$2.1$ in Claim~\ref{Claim4} 
(see Fig.~\ref{Fig27}(a)). 
This contradicts Lemma~\ref{PinWheelLemma}.

For Case~$2$. Since $\ell_6$ contains two white vertices $w_2, v^*_1$, 
the edge $e_4$ is not a terminal edge by Lemma~\ref{ConsecutiveTripletLemma4Chart}. 
But we cannot find any white vertex at which the edge $e_4$ is inward 
because the direction indicators 
for the arcs $\ell_6,\ell_8$ point to $D_4$ from $D_6$ 
and from $D_2$, respectively (see Fig.~\ref{Fig27}(b)).	
This is a contradiction.

Therefore, Claim~\ref{Claim5} holds.\hfill {$\square$}\vspace{1.5em} 

\begin{figure}[h]
\begin{center}
\includegraphics{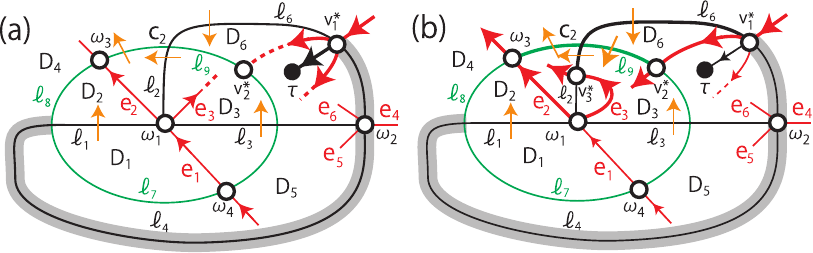}
\end{center}
\caption{\label{Fig28} 
Proof of Proposition~\ref{E1Terminal}. The terminal edge $\tau$ is an inner edge for $D_6$,
and is  outward at $v_1^*$.}
\end{figure}

{\it Proof of Proposition~\ref{E1Terminal}.}
Now we shall finish the proof of our proposition. 
By Claim~\ref{Claim4} and Claim~\ref{Claim5}, 
the terminal edge $\tau$ is an inner edge for $D_6$ and 
outward at the vertex $v_1^*$ on the arc $\ell_6$.
We start from Fig.~\ref{Fig28}(a).
By Lemma~\ref{BridgeLemma}, 
side edges of the terminal edge $\tau$ must 
contain a white vertex $ v^*_2$ on the arc $\ell_9$. 
Hence the direction indicator for the arc $\ell_9$ points to $D_3$ from $D_6$. 
Thus, the edge $e_3$ must contain a white vertex $ v^*_3$ 
on the arc $\ell_2$. 
Hence we can find a pinwheel with the vertices $w_3,v^*_1, v^*_2, v^*_3$ 
around the crossing $c_2$ (see Fig.~\ref{Fig28}(b)).
This contradicts Lemma~\ref{PinWheelLemma}.
Therefore, our proposition holds.
\hfill {$\square$}\vspace{1.5em}

\subsection{Case that the edge $e_1$ is a terminal edge}

\begin{figure}[h]
\begin{center}
\includegraphics{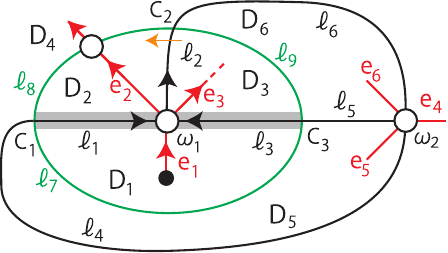}
\end{center}
\caption{\label{Fig29} 
The edge $e_1$ is a terminal edge.}
\end{figure}

Now we begin to prove our theorem. 
We start from the figure in Fig.~\ref{Fig29}.
Note that each of the arcs $\ell_1$ and $\ell_3$ 
does not contain any inner white vertices by Lemma~\ref{ConsecutiveTripletLemma4Chart}.

First of all, applying Proposition~\ref{E1Terminal} to $w_2$, 
one of the edges $e_4,e_5,e_6$ is a terminal edge. 
There are three cases.

{\bf Case A.} The edge $e_4$ is a terminal edge.

{\bf Case B.} The edge $e_5$ is a terminal edge.

{\bf Case C.} The edge $e_6$ is a terminal edge.

\begin{claim}
\label{Claim6} 
Neither Case A nor Case B occurs. 
\end{claim}

{\it Proof of Claim~\ref{Claim6}.}

\begin{figure}[h]
\begin{center}
\includegraphics{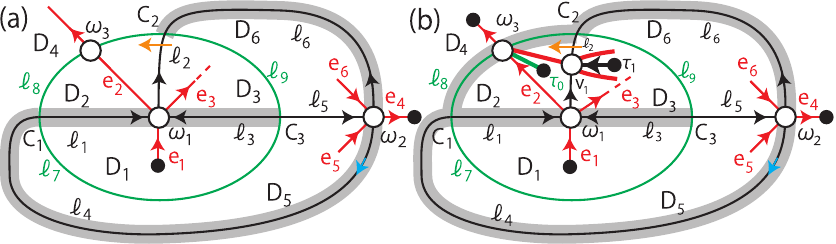}
\end{center}
\caption{\label{Fig30} 
Case A. The edge $e_4$ is a terminal edge.}
\end{figure}

For Case~A. 
Since $e_4$ is a terminal edge of label~$2$,
by Lemma~\ref{ConsecutiveTripletLemma4Chart} 
each of the arcs $\ell_4$ and $\ell_6$ 
contains no inner white vertices. 
Thus, the arc $\ell_4$ is outward at $w_2$, 
by considering $\ell_1\cup \ell_4$. 
Since the edge $e_4$ is a terminal edge, we have 
\begin{enumerate}
\item[(1)]
the edge $e_4$ is outward at $w_2$,  
the edge $e_5$ is inward at $w_2$, 
the edge $e_6$ is inward at $w_2$, 
the arc $\ell_5$ is inward at $w_2$, and 
the arc $\ell_6$ is outward at $w_2$ (see Fig.~\ref{Fig30}(a)).
\end{enumerate}

Since there does not exist any white vertex on $\ell_4\cup \ell_6$  except $w_2$ and 
since the edge $e_4$ is a terminal edge of label~$2$, 
the inner edge for $D_4$ containing the vertex $w_3$ 
is a terminal edge of label~$2$ outward at $w_3$. Further
\begin{enumerate}
\item[(2)] there does not exist any inner vertex of $\ell_8$ 
except $w_3$. 
\end{enumerate}

By Fact~\ref{Fact1},
there exists a terminal edge $\tau_0$ of label~$3$ containing $w_3$.
Then, the terminal edge $\tau_0$ is an inner edge for $D_2$, 
and the side edge of the terminal edge $\tau_0$ different from $e_2$ 
is outward at an inner vertex $v_1$ of the arc $\ell_2$.
Thus by Fact~\ref{Fact1},
there exists a terminal edge $\tau_1$ of label~$1$ containing $v_1$. 
The statement (2) implies that 
the terminal edge $\tau_1$ is an inner edge for the disk $D_3$ inward at the vertex $v_1$
(see Fig.~\ref{Fig30}(b)). 
By the same way as in the proof of Case~$2.1$ of Claim~\ref{Claim3} 
(see Fig.~\ref{Fig25}(c),(d)),
we can eliminate the white vertex $w_1$ by two C-I-M2 moves, 
and a C-III move.
This contradicts the minimality of the chart.
Hence Case~A does not occur.

\begin{figure}[h]
\begin{center}
\includegraphics{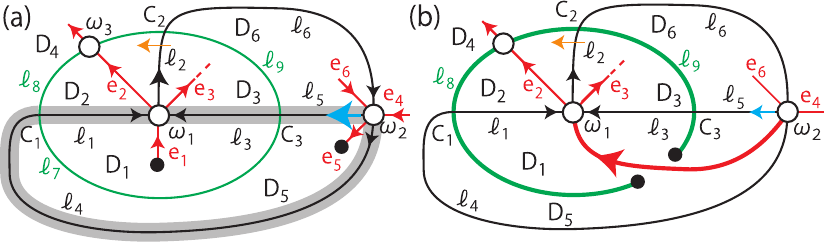}
\end{center}
\caption{\label{Fig31} 
(a) Case B. The edge $e_5$ is a terminal edge.
(b) The chart obtained by applying Cut Edge Lemma.}
\end{figure}

For Case B. 
By Lemma~\ref{ConsecutiveTripletLemma4Chart}, 
each of the arcs $\ell_1$ and $\ell_3$ 
contains no inner white vertices. 
Since $e_5$ is a terminal edge of label~$2$,
by Lemma~\ref{ConsecutiveTripletLemma4Chart} 
each of the arcs $\ell_4$ and $\ell_5$ 
contains no inner white vertices. 
Thus the arc $\ell_5$ is outward at $w_2$, 
by considering $\ell_3\cup \ell_5$. 
Since the edge $e_5$ is a terminal edge, we have 
\begin{enumerate}
\item[(3)]
 the edge $e_4$ is inward at $w_2$,  
the edge $e_5$ is outward at $w_2$, 
the edge $e_6$ is inward at $w_2$, 
the arc $\ell_4$ is outward at $w_2$, and 
the arc $\ell_6$ is inward at $w_2$ (see Fig.~\ref{Fig31}(a)).
\end{enumerate}

Now, there does not exist any inner vertex of $\ell_7$.
If not, then there exists an inner vertex of $\ell_7$.
By Fact~\ref{Fact1},
there exists a terminal edge $\tau$ of label $3$ containing the vertex.
By Lemma~\ref{BridgeLemma},
 the side edges of the terminal edge $\tau$
contain inner vertices of $\ell_1,\ell_3,\ell_4$ or $\ell_5$.
However, this contradicts
the fact that none of $\ell_1,\ell_3,\ell_4,\ell_5$
contain inner vertices.
Thus, there does not exist any inner vertex of $\ell_7$.

Looking at the orientation of the two terminal edges $e_1,e_5$, 
we can apply Cut Edge Lemma (Lemma~\ref{CutEdgeLemma})
to exchange the terminal edges
$e_1,e_5$ of label~$2$ to get  two new edges of label~$3$
containing black vertices 
(see Fig.~\ref{Fig31}(b)). 
By using the two edges of label~$3$, we eliminate the 
two crossings $c_1,c_3$ and the white vertices $w_1,w_3$ by two C-II moves and a C-I-M3 move. 
This contradicts the minimality of the chart.
Thus Case~B does not occur.

Therefore, Claim~\ref{Claim6} holds.\hfill {$\square$}\vspace{1.5em}

Now we investigate Case C.

\begin{figure}[h]
\begin{center}
\includegraphics{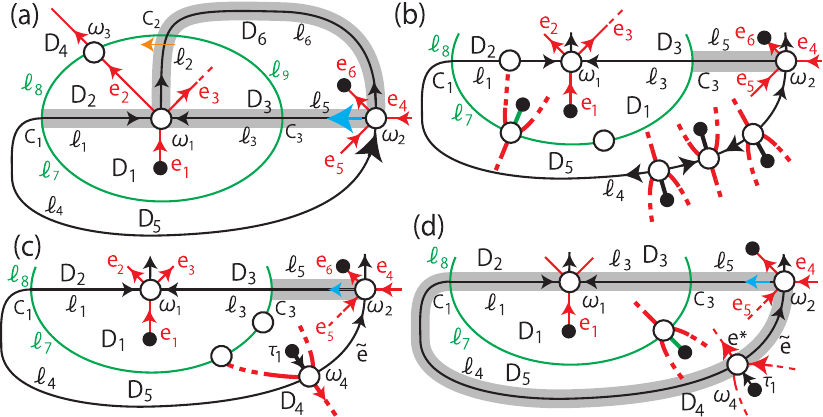}
\end{center}
\caption{\label{Fig32} 
Case C. (a) The edge $e_6$ is a terminal edge.
(b) Proof of Claim~\ref{Claim7} and Claim~\ref{Claim8}.
(c) Proof of Claim~\ref{Claim9}. The terminal edge $\tau_1$ is an inner edge for $D_5$.
(d) Proof of Claim~\ref{Claim10}. The terminal edge $\tau_1$ is an inner edge for $D_4$.
}
\end{figure}

For Case C. Since both of $e_1$ and $e_6$ are terminal edges of label~$2$,
by Lemma~\ref{ConsecutiveTripletLemma4Chart} 
\begin{enumerate}
\item[(a)]
each of the arcs $\ell_1$ and $\ell_3$ 
contains no inner white vertices, and
\item[(b)]
each of the arcs $\ell_5$ and $\ell_6$ 
contains no inner white vertices.
\end{enumerate}
Thus, the edge $\ell_5$ is outward at $w_2$, 
by considering $\ell_3\cup \ell_5$. 
Since the edge $e_6$ is a terminal edge, we have 
\begin{enumerate}
\item[(c)] the edge $e_4$ is inward at $w_2$,  
the edge $e_5$ is inward at $w_2$, 
the edge $e_6$ is outward at $w_2$, 
the arc $\ell_4$ is inward at $w_2$, and 
the arc $\ell_6$ is outward at $w_2$ (see Fig.~\ref{Fig32}(a)).
\end{enumerate}

\begin{claim}
\label{Claim7}
 There exists at most one inner vertex of the arc $\ell_7$.
\end{claim}

{\it Proof of Claim~\ref{Claim7}.}
If there exist at least two inner vertices of the arc $\ell_7$, 
then each inner vertex is contained in a terminal edge of label~$3$
by Fact~\ref{Fact1}.
Thus by Fact~\ref{Fact5}, 
there exists an inner terminal edge $\tau$ for $D_1$ of label~$3$.
Hence there exist two side edges of label~$2$ of the terminal edge $\tau$.
By Lemma~\ref{BridgeLemma},
one of the two side edges contains an inner vertex 
on $\ell_1\cup \ell_3$ different from $w_1$
 (see Fig.~\ref{Fig32}(b)). 
This contradicts Statement~(a) for Case C.
Hence, Claim~\ref{Claim7} holds.\hfill {$\square$}\vspace{1.5em}

\begin{claim}
\label{Claim8} 
There exists exactly one inner vertex of the arc $\ell_4$.
\end{claim}

{\it Proof of Claim~\ref{Claim8}.}
Since the arc $\ell_1$ is inward at $w_1$ and 
the arc $\ell_4$ is inward at $w_2$, 
the number of inner vertices on $\ell_1\cup \ell_4$ is odd, 
because the orientation on an IO-path changes 
at each white vertex 
(see Fig.~\ref{FigAnacanthousPath} and Fig.~\ref{Fig32}(b)). 
Since there does not exist any inner vertex on $\ell_1$
by Statement~(a) for Case~C, 
the number of inner vertices of the arc $\ell_4$ is odd.

Suppose that there exist at least three inner vertices of $\ell_4$.
By Fact~\ref{Fact5}, 
there exist at least four inner edges for $D_5$ of label~$2$ 
with inner vertices of $\ell_4$
(see Fig.~\ref{Fig32}(b)). 
By Fact~\ref{Fact4}, 
none of the four edges of label~$2$ are terminal edges.
Since the arc $\ell_4$ is an IO-path by Lemma~\ref{IO-pathLemma},
all of the four edges of label~$2$ are {\it inward} at the vertices on $\ell_4$ or 
all of the four edges of label~$2$ are {\it outward} at the vertices on $\ell_4$.
In either case,
there exist at least two inner vertices of $\ell_7$ 
by Statement~(b) for Case C. 
This contradicts Claim~\ref{Claim7}.
Hence, Claim~\ref{Claim8} holds.\hfill {$\square$}\vspace{1.5em}

Let $w_4$ be the inner vertex of the arc $\ell_4$.
By Fact~\ref{Fact1}, 
there exists a terminal edge $\tau_1$ of label~$1$ 
containing the vertex $w_4$.

\begin{claim}
\label{Claim9} 
The terminal edge $\tau_1$ is an inner edge for the disk $D_4$.
\end{claim}

{\it Proof of Claim~\ref{Claim9}.}
Suppose that the terminal edge $\tau_1$ is an inner edge 
for the disk $D_5$.
Let $\tilde{e}$ be the edge in the arc $\ell_4$ containing the vertices
$w_2$ and $w_4$.
Since the edge $\tilde{e}$ is inward at $w_2$ by Statement~(c) for Case~C, 
the edge $\tilde{e}$ is outward at $w_4$.
Thus, the terminal edge $\tau_1$ is inward at the vertex $w_4$. 
By Fact~\ref{Fact4} and Lemma~\ref{ConsecutiveTripletLemma4Chart},
there are three inner edges of label~$2$ for $D_5$ inward at vertices on
$\ell_4 \cup \ell_5$ none of which is a terminal edge.
Thus, there exist at least two inner vertices on the arc $\ell_7$ (see Fig.~\ref{Fig32}(c)).
This contradicts Claim~\ref{Claim7}.
Hence, Claim~\ref{Claim9} holds.\hfill {$\square$}\vspace{1.5em}

Let $e^*$ be the inner edge for $D_5$ containing the vertex $w_4$.
Then $e^*$ is of label $2$.
Since the arc $\ell_4$ is inward at $w_2$
by Statement~(c) for Case~C,
the terminal edge $\tau_1$ is inward at $w_4$  and 
$e^*$ is outward at $w_4$.

\begin{claim}
\label{Claim10}
 There does not exist any inner vertex of $\ell_7$.
\end{claim}

{\it Proof of Claim~\ref{Claim10}.}
Suppose there exists an inner vertex of $\ell_7$. 
By Fact~\ref{Fact1},
there exists a terminal edge $\tau$ of label $3$ containing the vertex. 
If the terminal edge $\tau$ is an inner edge for the disk $D_1$, 
then the side edges of $\tau$ contains 
an inner vertex of $\ell_1$ or $\ell_3$. 
This contradicts Statement~(a) for Case C.
If the terminal edge $\tau$ is an inner edge for the disk $D_5$,
then one of the side edges of $\tau$ 
contains a white vertex on
the arc $\ell_4\cup \ell_5$ different from $w_2$ and $w_4$ 
because the edge $e_5$ is inward at $w_2$ and 
the edge $e^*$ is outward at $w_4$ (see Fig.~\ref{Fig32}(d)). 
This contradicts Claim~\ref{Claim8} and Statement~(b) for Case~C. 
Hence, Claim~\ref{Claim10} holds.\hfill {$\square$}\vspace{1.5em}

\begin{figure}[h]
\begin{center}
\includegraphics{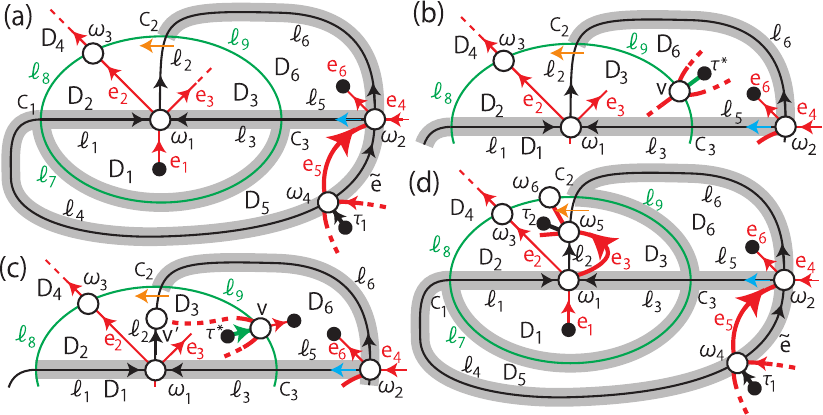}
\end{center}
\caption{\label{Fig33} 
(a) There exists exactly one inner vertex of the arc $\ell_4$, 
and
there does not exist any inner vertex of $\ell_7$.
(b) Case~$1$ of Claim~\ref{Claim11}. 
The terminal edge $\tau^*$ is an inner edge for $D_6$.
(c) Case~$2$ of Claim~\ref{Claim11}. 
The terminal edge $\tau^*$ is an inner edge for $D_3$.
(d) There does not exist any inner vertex of $\ell_9$.}
\end{figure}

Now we have a pseudo-chart in Fig.~\ref{Fig33}(a) 
by Claim~\ref{Claim10}.

\begin{claim}
\label{Claim11}
 There does not exist any inner vertex of the arc $\ell_9$.
\end{claim}

{\it Proof of Claim~\ref{Claim11}.}
Suppose that there exists an inner vertex $v$ of $\ell_9$.
By Fact~\ref{Fact1},
there exists a terminal edge $\tau^*$ of label 3 containing $v$.
There are two cases.

{\bf Case 1.} The terminal edge $\tau^*$ is an inner edge for the disk $D_6$.

{\bf Case 2.} The terminal edge $\tau^*$ is an inner edge for the disk $D_3$.

For Case~$1$. 
The terminal edge $\tau^*$ is an inner edge for the disk $D_6$.
One of the side edges of the terminal edge $\tau^*$ contains 
an inner vertex of $\ell_5$ or $\ell_6$
by Lemma~\ref{BridgeLemma}. 
This contradicts Statement~(b) for Case C 
because we are still in Case C (see Fig.~\ref{Fig33}(b)).

For Case~$2$. 
The terminal edge $\tau^*$ is an inner edge for the disk $D_3$. 
Let $e'$ be the inner edge for $D_6$ of label~$2$ containing $v$. 
If $e'$ is not a terminal edge,
then  $e'$ contains an inner vertex of $\ell_5$ or $\ell_6$,
because $\ell_9$ is an IO-path.
This contradicts Statement~(b) for Case C.
Hence, the edge $e'$ is a terminal edge of label~$2$.

If the terminal edge $\tau^*$ is outward at $v$, 
then the edge $e'$ is inward at $v$.
Then, we can get a new free edge by a C-I-M2 move between $e'$ and $e_6$.
This move increases the number of free edges without changing the number of crossings and white vertices, and hence decreases the $c$-complexity. 
This contradicts the minimality of the chart. 

Thus, the terminal edge $\tau^*$ is inward at $v$.
One of the side edges of $\tau^*$ contains an inner vertex $v'$ of 
$\ell_2$ or $\ell_3$. 
By Statement~(a) for Case C, 
the white vertex $v'$ is on $\ell_2$ 
(see Fig.~\ref{Fig33}(c)).
Since the terminal edge $\tau^*$ is inward at $v$, 
the side edge must be outward at $v'$. 
This is impossible because the direction indicator for $\ell_2$ 
points to $D_2$ from $D_3$.

Therefore, Claim~\ref{Claim11} holds.\hfill {$\square$}\vspace{1.5em}

By Claim~\ref{Claim11}, 
the edge $e_3$ contains an inner vertex $w_5$ of $\ell_2$.
By Fact~\ref{Fact1},
there exists a terminal edge $\tau_2$ of label $1$ containing $w_5$.
If the terminal edge $\tau_2$ is an inner edge for $D_3$,
then by the same way as in the proof of Case~$2.1$ of Claim~\ref{Claim3} 
(see Fig.~\ref{Fig25}(c),(d))
we can eliminate the white vertex $w_1$ by two C-I-M2 moves 
and a C-III move.
Hence, the terminal edge $\tau_2$ is an inner edge for $D_2$, 
and 
one of the side edges of the terminal edge $\tau_2$ contains 
an inner vertex $w_6$ of $\ell_8$ (see Fig.~\ref{Fig33}(d)).
Thus by Fact~\ref{Fact1},
there exists a terminal edge $\tau_3$ of label $3$ containing the white vertex $w_6$.

\begin{figure}[hbt]
\begin{center}
\includegraphics{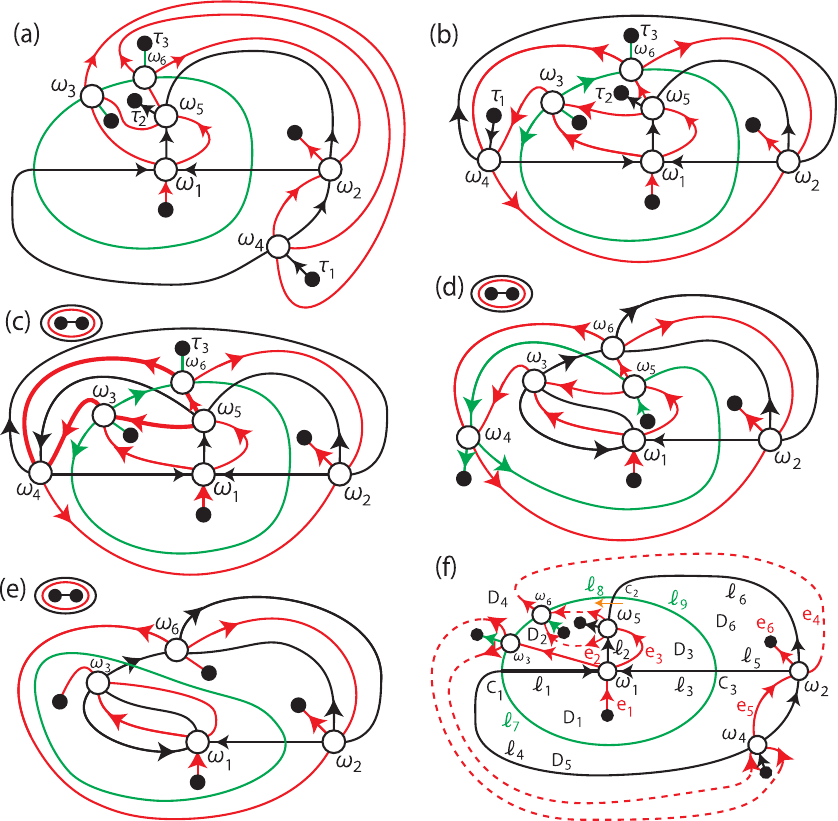}
\end{center}
\caption{\label{Fig34} 
(a),(b) $\tau_3$ is an inner edge for $D_4$.
(c),(d),(e) Charts which are C-move equivalent to the chart as shown (a).
(f) $\tau_3$ is an inner edge for $D_2$.}
\end{figure}

\begin{claim}
\label{Claim12}
 The terminal edge $\tau_3$ is an inner edge for $D_2$.
\end{claim}

{\it Proof of Claim~\ref{Claim12}.} 
Suppose that the terminal edge $\tau_3$ is an inner edge for $D_4$. 
Then by Fact~\ref{Fact5},
the terminal edge of label~$3$ at $w_3$ is an inner edge for $D_2$, and
we have the chart in Fig.~\ref{Fig34}(a)
which is isomorphic to the one in Fig.~\ref{Fig34}(b) 
on the sphere.
Move the terminal edge $\tau_1$ close to the terminal edge $\tau_2$. 
Apply a C-II move and a C-I-M2 move between the two terminal edges $\tau_1$ and $\tau_2$ to get one free edge and one crossing. 
Shift out the free edge from the square of label~$2$ 
with four vertices $w_6,w_4,w_3,w_5$ 
(see Fig.~\ref{Fig34}(c)).
Applying a C-I-M4 move, and two C-I-R2 moves 
(see Fig.~\ref{Fig34}(d)),
and then a C-III move to the terminal edge of label~$3$ at $w_4$.
Moreover, apply a C-III move to the terminal edge of label~$3$ at $w_5$,
we get the chart in Fig.~\ref{Fig34}(e) 
which contains two crossings and four white vertices 
(it actually becomes a ribbon chart at the end). 
Hence the chart in Fig.~\ref{Fig34}(a) 
is not a minimal chart.
Thus, Claim~\ref{Claim12} holds.\hfill {$\square$}\vspace{1.5em}

\begin{figure}[hbt]
\begin{center}
\includegraphics{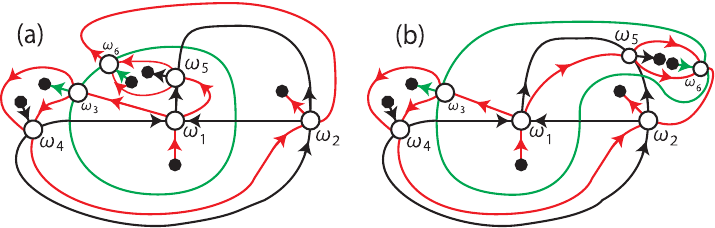}
\end{center}
\caption{\label{Fig37} 
4-charts describing a 2-twist-spun trefoil knot.}
\end{figure}

Therefore, the terminal edge $\tau_3$ is an inner edge for $D_2$.
Thus by Fact~\ref{Fact5},
the terminal edge of label~$3$ at $w_3$ is an inner edge for $D_4$.
We connect the corresponding ends in Fig.~\ref{Fig34}(f),
and move the white vertex $w_4$ close to the white vertex $w_3$ (see Fig.~\ref{Fig37}(a)).
Moreover, we move the white vertex $w_5$ close to the white vertex $w_2$ (see Fig.~\ref{Fig37}(b)).
Then we obtain the first 4-chart as shown in Fig.~\ref{Fig01}, 
up to the addition or deletion of free edges and hoops.
Hence the original chart is lor-equivalent to the chart describing the 2-twist-spun trefoil knot.
 This completes the proof of Theorem~\ref{MainTheorem}.

\section{Acknowledgement}

The second author is supported by JSPS KAKENHI Grant Number 24K06745.




\vspace{5mm}

\begin{minipage}{65mm}
{Teruo NAGASE
\\
{\small Tokai University \\
4-1-1 Kitakaname, Hiratsuka \\
Kanagawa, 259-1292 Japan\\
\\
nagase.teruo.t@tokai.ac.jp
}}
\end{minipage}
\begin{minipage}{65mm}
{Akiko SHIMA 
\\
{\small Department of Mathematics, 
\\
Tokai University
\\
4-1-1 Kitakaname, Hiratsuka \\
Kanagawa, 259-1292 Japan\\
shima-a@tokai.ac.jp
}}
\end{minipage}

\vspace{0.7cm}

{\bf List of terminology}\vspace{2mm}\\
{\small $
\begin{array}{ll||}
\text{$AB$-component} & p14 \\
\text{anacanthous path} & p8 \\
\text{anacanthous body $AB(\Gamma_1),AB(\Gamma_3)$} &  p14\\
\text{bigon} & p4 \\
\text{$c$-complexity $(c(\Gamma),w(\Gamma),-f(\Gamma),-b(\Gamma))$} & p4 \\
\text{chart} & p3 \\
\text{$c$-minimal chart} & p4 \\
\text{C-move equivalent} & p4 \\
\text{consecutive triplet} & p12 \\
\text{direction indicator} & p8 \\
\text{free edge} & p2 \\
\text{hoop} & p2 \\
\text{inner edge for a disk} & p8 \\
\text{inner vertex of an arc} & p8 \\
\text{inward} & p7 \\
\text{IO-path} & p18 \\
\end{array}
~~
\begin{array}{ll}
\text{linear} & p2 \\
\text{lor-equivalent} & p2 \\
\text{mal-cycle} & p8 \\
\text{middle arc} & p4 \\
\text{minimal chart} & p4 \\
\text{outward} & p7 \\
\text{oval nest} & p5 \\
\text{pinwheel} & p10 \\
\text{point at infinity $\infty$} & p6 \\
\text{pseudo-chart} & p15 \\
\text{ribbon chart} & p2 \\
\text{ribbon surface-link} &  p1\\
\text{side edge} & p12 \\
\text{simple hoop} & p5 \\
\text{terminal edge} & p5 \\
& \\
\end{array}
$}

\vspace{0.5cm}

{\bf List of notations}\vspace{2mm}\\
{\small $
\begin{array}{ll||}
\text{$\Gamma_m$} & p2 \\
\text{${\rm Cross}(\Gamma)$} & p2 \\
\text{$c(\Gamma)$} & p4\\
\text{$w(\Gamma)$} & p4 \\
\text{$f(\Gamma)$} & p4 \\
\end{array}
$~~
$\begin{array}{ll}
\text{$b(\Gamma)$} & p4\\
\text{${\rm Int}(X)$} & p7 \\
\text{$\partial X$} & p7 \\
\text{${\rm Cl}(X)$} & p7 \\
\text{$AB(\Gamma_1),AB(\Gamma_3)$}& p14\\
\end{array}
$
}

\end{document}